\documentclass[12pt,showkeys]{amsart}
\usepackage{color}
\usepackage{graphicx}
\usepackage{subfigure}
\usepackage{amssymb}
\usepackage{multirow}
\usepackage{amsmath}
\usepackage{multirow}
\usepackage{enumerate}
\usepackage{epstopdf}
\usepackage{cite,stmaryrd,txfonts}
\usepackage{tikz}

\input{undertilde}
\usetikzlibrary{snakes}

\usepackage{verbatim}

\usepackage{epic}

\usepackage{empheq}

\usepackage{ulem}

\newcommand{\bs}{\boldsymbol}

\def\vgm12{\bs{V}^{1+,2}_{\gamma,M}}

\newtheorem{theorem}{Theorem}

\newtheorem{lemma}[theorem]{Lemma}

\newcounter{mnote}
\let\oldmarginpar\marginpar
\renewcommand\marginpar[1]{\-\oldmarginpar[\raggedleft\footnotesize #1]%
  {\raggedright\footnotesize #1}}

\begin{document}

\title[A high accuracy nonconforming finite element scheme for TEP]{A high accuracy nonconforming finite element scheme for Helmholtz transmission eigenvalue problem}
\author{Yingxia Xi}
\author{Xia Ji}
\author{Shuo Zhang}
\address{School of Science, Nanjing University of Science and Technology, Nanjing 210094, China}
\email{xiyingxia@njust.edu.cn}
\address{LSEC, Institute of Computational Mathematics and Scientific/Engineering Computing, Academy of Mathematics and System Sciences, Chinese Academy of Sciences, Beijing 100190, People's Republic of China}
\thanks{The research of Y. Xi is supported in part by the National Natural Science Foundation of China with Grant No. 11901295, Natural Science Foundation of Jiangsu Province under BK20190431 and the Start-up Fund for Scientific Research, Nanjing University of Science and Technology (No. AE89991/109).}
\email{jixia@lsec.cc.ac.cn}
\address{LSEC, Institute of Computational Mathematics and Scientific/Engineering Computing, Academy of Mathematics and System Sciences, Chinese Academy of Sciences, Beijing 100190, People's Republic of China}
\thanks{The research of X. Ji is partially supported by the National Natural Science Foundation of China with Grant Nos. 11271018 and 91630313, and National Center for Mathematics and Interdisciplinary Sciences, Chinese Academy of Sciences.}
\email{szhang@lsec.cc.ac.cn}
\thanks{The research of S. Zhang is supported partially by the National Natural Science Foundation of China with Grant Nos. 11471026 and 11871465 and National Centre for Mathematics and Interdisciplinary Sciences, Chinese Academy of Sciences.}




\begin{abstract}
In this paper, we consider a cubic $H^2$ nonconforming finite element scheme $B_{h0}^3$ which does not correspond to a locally defined finite element with Ciarlet$'$s triple but admit a set of local basis
functions. For the first time, we deduce and write out the expression of basis functions explicitly.
Distinguished from the most nonconforming finite element methods, $(\delta\Delta_h\cdot,\Delta_h\cdot)$ with non-constant coefficient $\delta>0$ is coercive on the nonconforming  $B_{h0}^3$ space which makes it robust for numerical discretization. For fourth order eigenvalue problem, the $B_{h0}^3$ scheme can provide $\mathcal{O}(h^2)$ approximation for the eigenspace in energy norm and $\mathcal{O}(h^4)$ approximation for the eigenvalues.
We test the $B_{h0}^3$ scheme on the vary-coefficient bi-Laplace source and eigenvalue problem, further, transmission eigenvalue problem. Finally, numerical examples are presented to demonstrate the effectiveness of the proposed scheme.
\end{abstract}

\maketitle


\section{Introduction}
Recently the transmission eigenvalue problem has been attracting interests from many researchers. This problem arose in the inverse scattering theory for inhomogeneous medium and plays a key role in inverse scattering theory. The transmission eigenvalues can be used to obtain estimates for the physical characteristics of the hidden scatterer and have a variety of applications in inverse problem, such as target identification and nondestructive testing \cite{CakoniRen}\cite{HarrisCakoniSun}. Besides, transmission eigenvalues can also be used to design the invisible material \cite{JiLiu}.

Typically, for the scattering of time-harmonic acoustic waves by a bounded simply connected inhomogeneous medium $\Omega \subset \mathcal{R}^2$, the transmission eigenvalue problem is to find $k\in\mathcal{C}$, $\phi, \varphi\in H^2(\Omega)$ such that
\begin{equation*}
\left\{
\begin{array}{rcll}
\Delta \phi+k^2n(x)\phi&=&0,&{\rm in}\ \Omega,\\
\Delta \varphi+k^2\varphi&=&0, &{\rm in}\ \Omega,\\
\phi-\varphi&=&0,&{\rm on}\ \partial \Omega,\\
\frac{\partial \phi}{\partial \nu}-\frac{\partial \varphi}{\partial\nu}&=&0,&{\rm on}\ \partial \Omega,
\end{array}
\right.
\end{equation*}
where $n(x)$ is the index of refraction and $\nu$ is the unit outward normal to the boundary $\partial \Omega$. Typically, it's assumed that $n(x)>1$ or $0<n(x)<1$.

The transmission eigenvalue problem is non-self-adjoint and not covered by the standard theory of partial differential equations. It is numerically challenging because of the nonlinearity and the complicated spectral. Moreover, in most cases, the continuous problem degenerates with an infinite dimensional eigenspace associated with the zero eigenvalue, which has no physical meaning and makes it difficult to be solved. The first numerical study may be found in \cite{ColtonMonkSun2010} where three finite element methods were proposed. In \cite{Sun}, the author reformulates the transmission eigenvalue problem as the combination of a nonlinear function and a series of fourth order self-adjoint eigenvalue problems. The roots of the nonlinear function are the transmission eigenvalues., and an iterative method was proposed based on this. The rigorous convergence analysis was first given. But this method can only capture real eigenvalues.

To avoid the non-physical eigenspaces, introducing a new variable $u=\phi-\varphi\in H_0^2(\Omega)$, following the same procedure in \cite{JiSunTurner2012ACMTOM}, we can obtain the following fourth order equation
\begin{equation}\label{FourthEquation}
\big(\Delta +k^2n(x)\big)\frac{1}{n(x)-1}(\Delta +k^2)u=0.
\end{equation}
We remark that the above fourth order equation has eliminated the non-physical zero eigenvalue. Actually $k = 0$ implies $\displaystyle (\frac{1}{n(x)-1}\Delta u, \Delta u)=0$ and $u \in H_0^2(\Omega)$, and then we can obtain $u = 0$. The corresponding variational formulation of \eqref{FourthEquation} is to find $(k^2\neq 0,u)\in\mathbb{C}\times H_0^2(\Omega)$, such that
\begin{equation}\label{eq:vftep}
\left(\frac{1}{n(x)-1}(\Delta u+k^2u),\Delta v+k^2n(x)v\right)=0,\ \ \forall\, v\in H_0^2(\Omega).
\end{equation}
Let $\tau=k^2$ (we also call $\tau$ a transmission eigenvalue if $k$ is), the corresponding variational form is to find
$(\tau\neq 0,u)\in\mathbb{C}\times H_0^2(\Omega)$, such that
\begin{equation}\label{eq:vftep}
\left(\frac{1}{n(x)-1}(\Delta u+\tau u),\Delta v+\tau n(x)v\right)=0,\ \ \forall\, v\in H_0^2(\Omega).
\end{equation}
Here we consider the case $n(x)>1$ for illustration. For the case $0<n(x)<1$, it follows similarly.  Using Green formula, we can rewrite the original variational formulation (\ref{eq:vftep}) as
\begin{eqnarray}\label{Weak_Eigenvalue_A}
\mathcal{A}_{\tau}(u,v)&=&\tau \mathcal{B}(u,v),\ \ \ \forall v\in V,
\end{eqnarray}
where
\begin{eqnarray}
\mathcal{A}_{\tau}(u,v)&=&\Big(\frac{1}{n(x)-1}(\Delta u+\tau u),
(\Delta v+\tau v)\Big)+\tau^2\big(u,v\big),\label{A_Tau_Form}
\end{eqnarray}
and
\begin{eqnarray}
\mathcal{B}(u,v)&=&\big(\nabla u,\nabla v\big).\label{B_Form}
\end{eqnarray}
The bilinear form $\mathcal{A}_\tau(\cdot,\cdot)$ is coercive on $H_0^2(\Omega)\times H_0^2(\Omega)$, and the bilinear form $\mathcal{B}(\cdot,\cdot)$ is symmetric and nonnegative on $H_0^2(\Omega)\times H_0^2(\Omega)$\cite{CakoniHaddar,Sun}.

The finite element discretization of \eqref{FourthEquation} is natural. Many schemes, such as the Argyris element method \cite{ColtonMonkSun2010}, the (multi-level) BFS element method \cite{JiSunXie2014JSC}, the Morley element method \cite{JiXiXie2017,XiJi2018}, the modified Zienkiewciz element and the Morley-Zienkiewicz element \cite{YangHanBi2016} and other low complexity finite element methods including an interior penalty discontinuous Galerkin method using C$^0$ Lagrange elements (C$^0$IPG method) \cite{Geng2016},  and so on. There have also existed some mixed methods for this problem. The related works for mixed element method can be referred to \cite{Camano2018,ColtonMonkSun2010,JiSunTurner2012ACMTOM,XiJiZhang2018,YangBiLiHan2016}. The mixed scheme in \cite{Camano2018, JiSunTurner2012ACMTOM} which is close to the Ciarlet-Raviart discretization of biharmonic problem is based on Lagrange finite element method. For the nonzero transmission eigenvalues, this scheme is equivalent to the one proposed in \cite{ColtonMonkSun2010}. However, the scheme in \cite{Camano2018,JiSunTurner2012ACMTOM} can eliminate the zero transmission eigenvalue which has an infinite dimensional space and has no physical meaning. A mixed formulation in terms of three scaler fields and a spectral-mixed method are constructed in \cite{YangBiLiHan2016}. In \cite{XiJiZhang2018}, the authors propose a multi-level mixed formulation in terms of seven scaler fields. An equivalent linear mixed formulation of transmission eigenvalue problem which doesn't produce spurious modes even on non-convex domains is constructed. The proposed scheme admits a natural nested discretization, based on that a multi-level scheme is built. Optimal convergence rate and optimal computational cost can be obtained.

The finite element discretization of \eqref{Weak_Eigenvalue_A} looks immediate. While a $(\Delta\cdot,\Delta\cdot)$ bilinear form is used in the formulation, however, we have to note that $(\Delta_h\cdot,\Delta_h\cdot)$ is not coercive on general nonconforming finite element spaces. A standard approach is to enhance the bilinear from with $\alpha(\nabla^2\cdot,\nabla^2\cdot)$ for stabilisation, where $\alpha$ is a parameter. It is then not surprising that the choice of $\alpha$ may effect the performance of the scheme; a detailed illustration of the sensitivity of $\alpha$ can be found in Sections \ref{sec:morleybl} and  \ref{sec:morleyte}. To strengthen the robustness of the scheme, a finite element space which is of low degree and on which the bilinear form $(\Delta_h\cdot,\Delta_h\cdot)$ is coercive is needed.

In this paper, we introduce a new scheme for the Helmholtz transmission eigenvalue problem. Basically, we adapt onto \eqref{Weak_Eigenvalue_A} a piecewise cubic finite element space $B^3_{h0}$ introduced in \cite{Zhang2018,Zhang2019}. It is proved that $B^3_{h0}$ provide $\mathcal{O}(h^2)$ accuracy on both approximation error in broken $H^2$ error and consistency error associated to the biharmonic operator. Moreover, it is proved in \cite{Zhang2019} that
$$(\Delta_hu_h,\Delta_hv_h)=(\nabla_h^2u_h,\nabla_h^2v_h),\ \forall\, u_h,v_h\in B^3_{h0}.$$
Thus a finite element scheme based on $B^3_{h0}$ for the transmission eigenvalue problem can provide $\mathcal{O}(h^2)$ approximation for the eigenspace in energy norm and $\mathcal{O}(h^4)$ approximation for the eigenvalues. Numerical experiments of this paper verify this.

The space $B^3_{h0}$ does not correspond to a finite element defined by Ciarlet's triple, however, it admits a set of local basis functions\cite{Zhang2018}. By following the procedure given in \cite{Zhang2018,Zhang2019}, the finite element scheme designed in this paper can be implemented without knowing the basis functions of $B^3_{h0}$. However, in the order that basic algorithms can be used, the local basis functions are still in need, and we figure out them in this paper.

The rest of this paper is organized as follows. In Section \ref{ref:biLap}, we study the finite element space $B^3_{h0}$ and its utilization for the bi-Laplacian operator. We particularly figure out its local basis functions and illustrate the performance of the scheme with numerical examples. An illustration about the Morley element onto the model problem is also given for comparison. Section \ref{sec:te} is devoted to the Helmholtz transmission eigenvalue problem. Numerical experiments are given, including those of the Morley element for comparison. Finally, some concluding remarks are given in Section \ref{sec:conc}.

\section{A high-accuracy scheme for bi-Laplacian problem with varying coefficient}
\label{ref:biLap}

In this section, we first consider the following fourth order eigenvalue problem
\begin{equation}\label{Non-constantBiharEigenProblem}
\left\{
\begin{array}{rcl}
\Delta (\delta \Delta u)&=&\lambda u,\ \ \ \ \rm{in}\ \ \Omega,\\
u&=&0,\ \ \ \ \rm{on}\ \ \partial\Omega,\\
\frac{\partial u}{\partial n}&=&0,\ \ \ \ \rm{on}\ \ \partial\Omega,
\end{array}
\right.
\end{equation}
where $\delta(x)$ is a bounded smooth non-constant function and $\delta\geqslant \delta_{\min}>0$.

\subsection{A piecewise cubic finite element space and its structure}


Before introducing this finite element, we introduce some notations.
We assume $\mathcal{T}_h$ a shape regular mesh over $\Omega$ with mesh size $h$. Denote $\mathcal{X}_h$, $\mathcal{X}_h^i$, $\mathcal{X}_h^b$, $\mathcal{E}_h$, $\mathcal{E}_h^i$, $\mathcal{E}_h^b$ the vertices, interior vertices, boundary vertices, the set of edges, interior edges and boundary edges, respectively. For any edge $e\in \mathcal{E}_h$, denote the unit normal vector of $e$ by $\bf{n}_e$.
For a fixed element $T\in\mathcal{T}_h$, we denote $\mathcal{P}_k(T)$ the polynomial space of degree less than or equal to k and $|T|$ means the area measurement of element $T$. On an edge $e$, $\mathcal{P}_k(e)$ and $|e|$ are defined similarly.
The barycentre coordinates are denoted as usual by $\lambda_i(i=1,2,3).$

The nonconforming finite element space $B_h^3$ can be defined
as follows: (\cite{Zhang2019,Zhang2018})
\begin{eqnarray*} 
B_h^3&=&\big\{v\in L^2(\Omega)\mid\ v|_T\in \mathcal{P}_3(T),\ v\ {\rm is\ continuous\ at\ vertices}\ a\in\mathcal{X}_h  \rm{\ and}\nonumber\\
&&\ \int_{e}\llbracket v\rrbracket\ ds =0,\ {\rm  and}\ \int_{e}p_e\llbracket \partial_n v \rrbracket\ ds =0,\ \forall\ p_e\in P_1(e),\ \forall\ e\in \mathcal{E}_h^i,\ \forall\ T\in\mathcal{T}_h\big\}
\end{eqnarray*}
where $\llbracket v\rrbracket$ represents the jump of the scalar function $v$ across e,
and
\begin{eqnarray*} 
B_{h0}^3&=&\big\{v\in B_h^3\mid\ v(a)=0, a\in \mathcal{X}_h^{b};\ \int_{e} v\ ds =0,\ {\rm  and}\ \int_{e}p_e \partial_n v \ ds =0,\ \forall\ p_e\in P_1(e),\ \forall\ e\in \mathcal{E}_h^b\}.
\end{eqnarray*}

\begin{lemma}\cite{Zhang2018,Zhang2019}
$\inf_{w_h\in B^3_{h0}}|w-w_h|_{2,h}\leqslant Ch^k|w|_{2+k,\Omega},\ \forall\ w\in H^2_0(\Omega)\cap H^{k+2}(\Omega),\ \ k=1,2$.
\end{lemma}


\paragraph{\bf Local basis functions of $B^3_{h0}$}

The space $B_{h0}^3$ does not correspond to a locally defined finite element with Ciarlet$'$s triple. However, it is pointed out that the space admits a set of local basis functions. In the following, we will deduce and write out the expressions of basis functions in detail.
The derivation is based on the thought raised in \cite{Zhang2018} and we need the following results.
\begin{lemma}(\cite{Zhang2018})
$B_{h0}^3$ admits a set of basis functions with vertex-patch-based supports.
\end{lemma}

The following lemma involves the vector-valued finite element spaces $\widetilde{S}_{h0}^2(rot,w_0)$, $\widetilde{G}_{h0}^2(rot_h,0)$ of which the definition concerns a series of definitions of associated finite element spaces. It's omitted here and the author can refer the detail to \cite{Zhang2018}. And we use "$~\widetilde{}~$" for vector valued quantities in the following. And
$\widetilde{\varphi}^1$, $\widetilde{\varphi}^2$ are the two components of the quantity $\widetilde{\varphi}$.

\begin{lemma}(\cite{Zhang2018})\label{DeduceBh03}
Define an operator $\mathcal{F}_h: \widetilde{S}_{h0}^2(rot,w_0)\longrightarrow \widetilde{G}_{h0}^2(rot_h,0)$ by
\begin{equation}
\mathcal{F}_h \widetilde{\varphi}_h=\widetilde{\varphi}_h+\widetilde{\phi}_h,\ \ \forall\widetilde{\varphi}_h\in\widetilde{S}_{h0}^2(rot,w_0),\ \ \widetilde{\phi}_h\in \widetilde{B}_{h0}^2,\ s.t.\ \ rot_h(\mathcal{F}_h\widetilde{\varphi}_h)=0,  \nonumber
\end{equation}
where $\widetilde{B}_{h0}^2=\{\widetilde{\phi}_h:\ (\widetilde{\phi}_h|_T)^j\in \ span\{(\lambda_1^2+\lambda_2^2+\lambda_3^2)-2/3\},\ j=1,2,\ \forall\ T\in\mathcal{T}_h\}$.
And define $(\nabla^{-1})_h: \widetilde{G}_{h0}^2(rot_h,0)\longrightarrow B_{h0}^3$, then $(\nabla^{-1})_h\circ\mathcal{F}_h: \widetilde{S}_{h0}^2(rot,w_0)\longrightarrow B_{h0}^3$ is bijective and preserves support.
\end{lemma}

From Lemma \ref{DeduceBh03}, it can be observed that there are three steps in the derivations of basis functions. We orderly construct the basis functions in $\widetilde{S}_{h0}^2(rot,w_0)$, $\widetilde{G}_{h0}^2(rot_h,0)$ and $B_{h0}^3$.
Before introducing the derivation, we give some definitions. For $a\in\mathcal{X}_h$, denote by $P_a$ the union of triangles of which $a$ is a vertex, namely the patch associated with $a$; for $e\in\mathcal{E}_h$, denote by $P_e$ the patch associated with $e$.

\textbf{First, we consider constructing the basis functions in $\widetilde{S}_{h0}^2(rot,w_0)$ with vertex-patch-based supports.} On every vertex (e.g: denoted by $a$), three basis functions are associated, which are labelled as $\widetilde{\varphi}_a^x,\ \widetilde{\varphi}_a^y,\ \widetilde{\varphi}_{P_a}$. And on every edge (e.g: denoted by $e$), one basis function is associated, which is labelled as $\widetilde{\varphi}_{e}$. For every basis function associated with an interior vertex $a$, its restriction on a cell $T$ such that $a$ is a node of $T$. For every basis function associated with an edge, its restriction on a cell $T$ such that $e$ is an edge of $T$. Then we can only focus on an element and give out the basis functions.

For the construction of $\widetilde{\varphi}_a^x,\ \widetilde{\varphi}_a^y,\ \widetilde{\varphi}_{P_a}$ and $\widetilde{\varphi}_{e}$, we follow the thought in \cite{Zhang2018} and have the guaranteed theoretical result.
\begin{lemma}(\cite{Zhang2018})
The set $\{\widetilde{\varphi}_a^x,\ \widetilde{\varphi}_a^y,\ \widetilde{\varphi}_{P_a},\ \widetilde{\varphi}_{e}\}_{a\in \mathcal{X}_h^i,\ e\in\mathcal{E}_h^i}$ forms a basis of $\widetilde{S}_{h0}^2(rot,w_0)$.
\end{lemma}

For a fixed element $T\in\mathcal{T}_h$, the vertex denoted by $i(i=1,2,3)$, the opposite side of vertex $i$ denoted by $e_i$.
For the vertex $i$ and its opposite edge $e_i$, the associated basis functions are as follows.
\begin{equation}
\bullet~\widetilde{\varphi}_i^x=((\widetilde{\varphi}_i^x)^1,(\widetilde{\varphi}_i^x)^2)^T=(\lambda_i-3\lambda_i\lambda_j-3\lambda_i\lambda_k,0)^T=(3\lambda_j^2+6\lambda_j\lambda_k+3\lambda_k^2-4\lambda_j-4\lambda_k+1,0)^T, \nonumber
\end{equation}
\begin{equation}
\bullet~\widetilde{\varphi}_i^y=((\widetilde{\varphi}_i^y)^1,(\widetilde{\varphi}_i^y)^2)^T=(0,\lambda_i-3\lambda_i\lambda_j-3\lambda_i\lambda_k)^T=(0,3\lambda_j^2+6\lambda_j\lambda_k+3\lambda_k^2-4\lambda_j-4\lambda_k+1)^T,
\nonumber
\end{equation}
\begin{equation}
\bullet~\widetilde{\varphi}_e^i=((\widetilde{\varphi}_e^i)^1,(\widetilde{\varphi}_e^i)^2)^T=\frac{6\lambda_j\lambda_k}{|e_i|}(-\tau_2(e_i),\tau_1(e_i))^T,
\nonumber
\end{equation}
\begin{multline}
\bullet~\widetilde{\varphi}_{P_i}=((\widetilde{\varphi}_{P_i})^1,(\widetilde{\varphi}_{P_i})^2)^T=\frac{6\lambda_i\lambda_j}{|e_k|}(\tau_1(e_k),\tau_2(e_k))^T+\frac{6\lambda_i\lambda_k}{|e_j|}(\tau_1(e_j),\tau_2(e_j))^T
\\=\frac{6(1-\lambda_j-\lambda_k)\lambda_j}{|e_k|}(\tau_1(e_k),\tau_2(e_k))^T+\frac{6(1-\lambda_j-\lambda_k)\lambda_k}{|e_j|}(\tau_1(e_j),\tau_2(e_j))^T,
\nonumber
\end{multline}
where $i,j,k$ satisfy the cyclic coordinate.

\textbf{Second, we consider constructing the basis functions in $\widetilde{G}_{h0}^2(rot_h,0)$.} By Lemma \ref{DeduceBh03} and its process of proof in \cite{Zhang2018}, it's easy to verify the following conclusion.
\begin{lemma}
Under the assumption that $\{\widetilde{\varphi}_a^x,\ \widetilde{\varphi}_a^y,\ \widetilde{\varphi}_{P_a},\ \widetilde{\varphi}_{e}\}_{a\in \mathcal{X}_h^i,\ e\in\mathcal{E}_h^i}$ forms a basis of $\widetilde{S}_{h0}^2(rot,w_0)$, then $\{\mathcal{F}_h\widetilde{\varphi}_a^x,\ \mathcal{F}_h\widetilde{\varphi}_a^y,\ \mathcal{F}_h\widetilde{\varphi}_{P_a},\ \mathcal{F}_h\widetilde{\varphi}_{e}\}_{a\in \mathcal{X}_h^i,\ e\in\mathcal{E}_h^i}$ forms a basis of $\widetilde{G}_{h0}^2(rot,0)$ and
\begin{equation}
supp(\mathcal{F}_h\widetilde{\varphi}_a^x)\subset supp(\widetilde{\varphi}_a^x),\  supp(\mathcal{F}_h\widetilde{\varphi}_a^y)\subset supp(\widetilde{\varphi}_a^y),\
supp(\mathcal{F}_h\widetilde{\varphi}_{P_a})\subset supp(\widetilde{\varphi}_{P_a}),\
supp(\mathcal{F}_h\widetilde{\varphi}_{e})\subset supp(\widetilde{\varphi}_{e}). \nonumber
\end{equation}
\end{lemma}
The above lemma tells us that $\mathcal{F}_h$ can preserve the linear independence and the support of basis functions. Then, for an element $T\in\mathcal{T}_h$, $\mathcal{F}_h\widetilde{\varphi}_i^x,\ \mathcal{F}_h\widetilde{\varphi}_i^y,\ \mathcal{F}_h\widetilde{\varphi}_{P_i},\ \mathcal{F}_h\widetilde{\varphi}_{e_i}(i=1,2,3)$ are the corresponding basis functions in $\widetilde{G}_{h0}^2(rot_h,0)$.

Denote $\phi_T=\lambda_1^2+\lambda_2^2+\lambda_3^2-2/3$ and $\widehat{\varphi}_i^x\triangleq\mathcal{F}_h \widetilde{\varphi}_i^x\in\widetilde{G}_{h0}^2(rot_h,0)$. By Lemma \ref{DeduceBh03}, we assume
\begin{equation}
\widehat{\varphi}_i^x=\widetilde{\varphi}_i^x+(\alpha_i^x,\beta_i^x)^T\phi_T,\ \ \ \ rot_h \widehat{\varphi}_i^x=0. \nonumber
\end{equation}

By calculation, we can obtain
\begin{equation}
\alpha_i^x=\frac{[(\lambda_k)_x-(\lambda_j)_x][(\lambda_k)_y+(\lambda_j)_y]}{(\lambda_k)_y(\lambda_j)_x-(\lambda_k)_x(\lambda_j)_y},\ \ \ \ \beta_i^x=\frac{[(\lambda_k)_y-(\lambda_j)_y][(\lambda_k)_y+(\lambda_j)_y]}{(\lambda_k)_y(\lambda_j)_x-(\lambda_k)_x(\lambda_j)_y}. \nonumber
\end{equation}
Similarly, for $\mathcal{F}_h\widetilde{\varphi}_i^y=\widetilde{\varphi}_i^y+(\alpha_i^y,\beta_i^y)^T\phi_T$, $\mathcal{F}_h\widetilde{\varphi}_e^i=\widetilde{\varphi}_e^i+(\alpha_e^i,\beta_e^i)^T\phi_T$,
$\mathcal{F}_h\widetilde{\varphi}_{P_i}=\widetilde{\varphi}_{P_i}+(\alpha_{P_i},\beta_{P_i})^T\phi_T$, we have
\begin{equation}
\alpha_i^y=\frac{[(\lambda_j)_x-(\lambda_k)_x][(\lambda_j)_x+(\lambda_k)_x]}{(\lambda_k)_y(\lambda_j)_x-(\lambda_k)_x(\lambda_j)_y},\ \ \ \ \beta_i^y=\frac{[(\lambda_j)_y-(\lambda_k)_y][(\lambda_j)_x+(\lambda_k)_x]}{(\lambda_k)_y(\lambda_j)_x-(\lambda_k)_x(\lambda_j)_y}. \nonumber
\end{equation}
\begin{equation}
\alpha_e^i=(\lambda_i)_x\frac{6|T|}{|e_i|^2},\ \ \ \
\beta_e^i=(\lambda_i)_y\frac{6|T|}{|e_i|^2}. \nonumber
\end{equation}
\begin{equation}
\alpha_{P_i}=\frac{-12|T|\nabla \lambda_j\nabla \lambda_k\left\{(\frac{(\lambda_k)_x}{|e_k|^2}+\frac{(\lambda_j)_x}{|e_j|^2})\right\}-\frac{3}{|T|}\left\{(\lambda_j)_x+(\lambda_k)_x\right\}}
{2(\lambda_k)_y(\lambda_j)_x-2(\lambda_k)_x(\lambda_j)_y}, \nonumber
\end{equation}
\begin{equation}
\beta_{P_i}=\frac{-12|T|\nabla \lambda_j\nabla \lambda_k\left\{(\frac{(\lambda_k)_y}{|e_k|^2}+\frac{(\lambda_j)_y}{|e_j|^2})\right\}-\frac{3}{|T|}\left\{(\lambda_j)_y+(\lambda_k)_y\right\}}
{2(\lambda_k)_y(\lambda_j)_x-2(\lambda_k)_x(\lambda_j)_y}. \nonumber
\end{equation}

\textbf{Here, we consider constructing the basis functions in $B_{h0}^3$.} By Lemma \ref{DeduceBh03} and its process of proof in \cite{Zhang2018}, it's easy to verify the following conclusion.
\begin{lemma}
Under the assumption that $\{\widetilde{\varphi}_a^x,\ \widetilde{\varphi}_a^y,\ \widetilde{\varphi}_{P_a},\ \widetilde{\varphi}_{e}\}_{a\in \mathcal{X}_h^i,\ e\in\mathcal{E}_h^i}$ forms a basis of $\widetilde{S}_{h0}^2(rot,w_0)$, then $\{(\nabla^{-1})_h\circ\mathcal{F}_h\widetilde{\varphi}_a^x,\ (\nabla^{-1})_h\circ\mathcal{F}_h\widetilde{\varphi}_a^y,\
(\nabla^{-1})_h\circ\mathcal{F}_h\widetilde{\varphi}_{P_a},\ (\nabla^{-1})_h\circ\mathcal{F}_h\widetilde{\varphi}_{e}\}_{a\in \mathcal{X}_h^i,\ e\in\mathcal{E}_h^i}$ forms a basis of $B_{h0}^3$ and
\begin{eqnarray}\nonumber
&supp\left((\nabla^{-1})_h\circ\mathcal{F}_h\widetilde{\varphi}_a^x\right)\subset supp(\widetilde{\varphi}_a^x),\  &supp\left((\nabla^{-1})_h\circ\mathcal{F}_h\widetilde{\varphi}_a^y\right)\subset supp(\widetilde{\varphi}_a^y),\
\\
&supp\left((\nabla^{-1})_h\circ\mathcal{F}_h\widetilde{\varphi}_{P_a}\right)\subset supp(\widetilde{\varphi}_{P_a}),\
&supp\left((\nabla^{-1})_h\circ\mathcal{F}_h\widetilde{\varphi}_{e}\right)\subset supp(\widetilde{\varphi}_{e}). \nonumber
\end{eqnarray}
\end{lemma}

Denote $w_i^x=(\nabla^{-1})_h\circ\mathcal{F}_h\widetilde{\varphi}_i^x$, $w_i^y=(\nabla^{-1})_h\circ\mathcal{F}_h\widetilde{\varphi}_i^y$, $w_e^i=(\nabla^{-1})_h\circ\mathcal{F}_h\widetilde{\varphi}_e^i$, $w_{P_i}=(\nabla^{-1})_h\circ\mathcal{F}_h\widetilde{\varphi}_{P_i}$. By calculation, the corresponding basis functions in $B_{h0}^3$ are as follows.\begin{multline}
w_i^x(\lambda_j,\lambda_k)=
-\xi_k\left\{(\frac{\lambda_j^3}{3}-\lambda_j^2+\frac{2}{3}\lambda_j)+(\frac{2}{3}\lambda_k^3-\lambda_k^2+\frac{\lambda_k}{3})+(2\lambda_j^2\lambda_k+\lambda_j\lambda_k^2-2\lambda_j\lambda_k)\right\}+
\\
\xi_j\left\{(\frac{2}{3}\lambda_j^3-\lambda_j^2+\frac{\lambda_j}{3})+(\frac{\lambda_k^3}{3}-\lambda_k^2+\frac{2}{3}\lambda_k)+(2\lambda_j\lambda_k^2+\lambda_j^2 \lambda_k-2\lambda_j\lambda_k)\right\},   \nonumber
\end{multline}
\begin{multline}
w_i^y(\lambda_j,\lambda_k)=-\eta_k\left\{(\frac{\lambda_j^3}{3}-\lambda_j^2+\frac{2}{3}\lambda_j)+(\frac{2}{3}\lambda_k^3-\lambda_k^2+\frac{\lambda_k}{3})+(2\lambda_j^2\lambda_k+\lambda_j\lambda_k^2-2\lambda_j\lambda_k)\right\}+
\\
\eta_j\left\{(\frac{2}{3}\lambda_j^3-\lambda_j^2+\frac{\lambda_j}{3})+(\frac{\lambda_k^3}{3}-\lambda_k^2+\frac{2}{3}\lambda_k)+(2\lambda_j\lambda_k^2+\lambda_j^2 \lambda_k-2\lambda_j\lambda_k)\right\}, \nonumber
\end{multline}
\begin{multline}
w_e^i(\lambda_j,\lambda_k)=-\frac{6|T|}{|e_i|^2}\left\{(\frac{2}{3}\lambda_j^3-\lambda_j^2+\frac{\lambda_j}{3})+(\frac{2}{3}\lambda_k^3-\lambda_k^2+\frac{\lambda_k}{3})+(2\lambda_j^2\lambda_k+2\lambda_j\lambda_k^2-2\lambda_j\lambda_k) \right\}, \nonumber
\end{multline}
\begin{multline}
w_{P_i}(\lambda_j,\lambda_k)=-3(\eta_j\eta_k+\xi_j\xi_k)\left\{\frac{1}{|e_j|^2}(\frac{2}{3}\lambda_j^3-\lambda_j^2+\frac{\lambda_j}{3})+\frac{1}{|e_k|^2}(\frac{2}{3}\lambda_k^3-\lambda_k^2+\frac{\lambda_k}{3})+\right\}
\\
+6\left\{(-\frac{2}{3}\lambda_j^3+\lambda_j^2-\frac{\lambda_j}{6})+(-\frac{2}{3}\lambda_k^3+\lambda_k^2-\frac{\lambda_k}{6})+  (\lambda_j\lambda_k-\lambda_j^2\lambda_k-\lambda_j\lambda_k^2)\right\}-1, \nonumber
\end{multline}
where $i=1, 2, 3$ which correspond to three vertices of a triangular element and $\xi_i=x_j-x_k$, $\eta_i=y_j-y_k$, $i, j, k$ satisfy the cyclic coordinate.

\subsection{A second order computational scheme for bi-Laplacian source and eigenvalue problems}

The bi-Laplacian source problem is to find $u$ satisfying
\begin{equation}\label{BiharBoundaryProblem}
\left\{
\begin{array}{rcl}
\Delta(\delta\Delta u)&=&f,\ \ \ \ \rm{in}\ \ \Omega,\\
u&=&0,\ \ \ \ \rm{on}\ \ \partial\Omega,\\
\frac{\partial u}{\partial n}&=&0,\ \ \ \ \rm{on}\ \ \partial\Omega.
\end{array}
\right.
\end{equation}
A finite element scheme for \eqref{BiharBoundaryProblem} is defined as: find $u_h\in B^3_{h0}$, such that
\begin{equation}\label{Non-constantBiharEigenProblemdis}
(\delta\Delta_hu_h,\Delta_hv_h)=(f,v_h),\quad\forall\,v_h\in B^3_{h0}.
\end{equation}
\begin{theorem}
Let $u\in H^4(\Omega)\cap H^2_0(\Omega)$ be the solution of \eqref{BiharBoundaryProblem}, and $u_h$ be the solution of \eqref{Non-constantBiharEigenProblemdis}, respectively. Then
$$
|u-u_h|_{2,h}\leqslant Ch^{k}|u|_{2+k,\Omega},\ \ \ k=1,2,
$$
and
$$
|u-u_h|_{1,h}\leqslant Ch^3|u|_{4,\Omega},\ \ \mbox{when}\ \Omega\ \mbox{is\ convex}.
$$
\end{theorem}

The finite element space $B^3_{h0}$ leads immediately to a high-accuracy scheme for the eigenvalue problem of bi-Laplacian equation.

\subsection{Numerical experiments}

\subsubsection{For source problems}

\textbf{Example 1.} Consider the bi-Laplacian source problem \eqref{BiharBoundaryProblem} with constant coefficient $\delta=1$ on square domain $\Omega_1=[0,1]^2$ with
$$f=-4\pi^4\left(\cos(2\pi x)+\cos(2\pi y)-4\cos(2\pi x)\cos(2\pi y)\right).$$
The exact solution is $u(x,y)=\sin(\pi x)^2\sin(\pi y)^2$.

\textbf{Example 2.} Consider the bi-Laplacian source problem \eqref{BiharBoundaryProblem} with constant coefficient $\delta=1$ on triangle domain $\Omega_2$ whose vertices are given by $(0,0), (1,0), (0,1)$.
And we consider $f=72(x+y)^2-48(x+y)+8$ for which the exact solution is $u(x,y)=x^2y^2(1-x-y)^2$.

Here we test on \textbf{Example 1} and \textbf{Example 2}, respectively. The mesh size of the
initial mesh is $h_0= \frac{1}{2}$. Six levels of uniformly refined triangular meshes are generated for numerical experiments and $h_k= h_{k-1}/2,\ k = 1,2,3,4,5,6$. The finest degrees of freedom (short for DOFs) for \textbf{Example 1} are 97283. The refinest DOFs for \textbf{Example 2} are 48387. We discretize by the second order computational scheme corresponding to $B_{h0}^3$ space. For each series of meshes, we obtain the numerical solution $u_{h_k}$. The convergent orders measured by $h_2,\ h_1,\ L_2$ norms respectively are computed by
\begin{equation}
log_2\left(\frac{\|u-u_{h_{k}}\|_{h_2}}{\|u-u_{h_{k-1}}\|_{h_2}}\right),\ \ \ \,k=2,3,4,5,6,\nonumber
\end{equation}
\begin{equation}
log_2\left(\frac{\|u-u_{h_{k}}\|_{h_1}}{\|u-u_{h_{k-1}}\|_{h_1}}\right),\ \ \ \,k=2,3,4,5,6,\nonumber
\end{equation}
and
\begin{equation}
log_2\left(\frac{\|u-u_{h_{k}}\|_{L_2}}{\|u-u_{h_{k-1}}\|_{L_2}}\right),\ \ \ \,k=2,3,4,5,6.\nonumber
\end{equation}

\begin{figure}[ht]
\centering
\subfigure[]{\label{SourceProblem_square}\includegraphics[scale=0.4503]{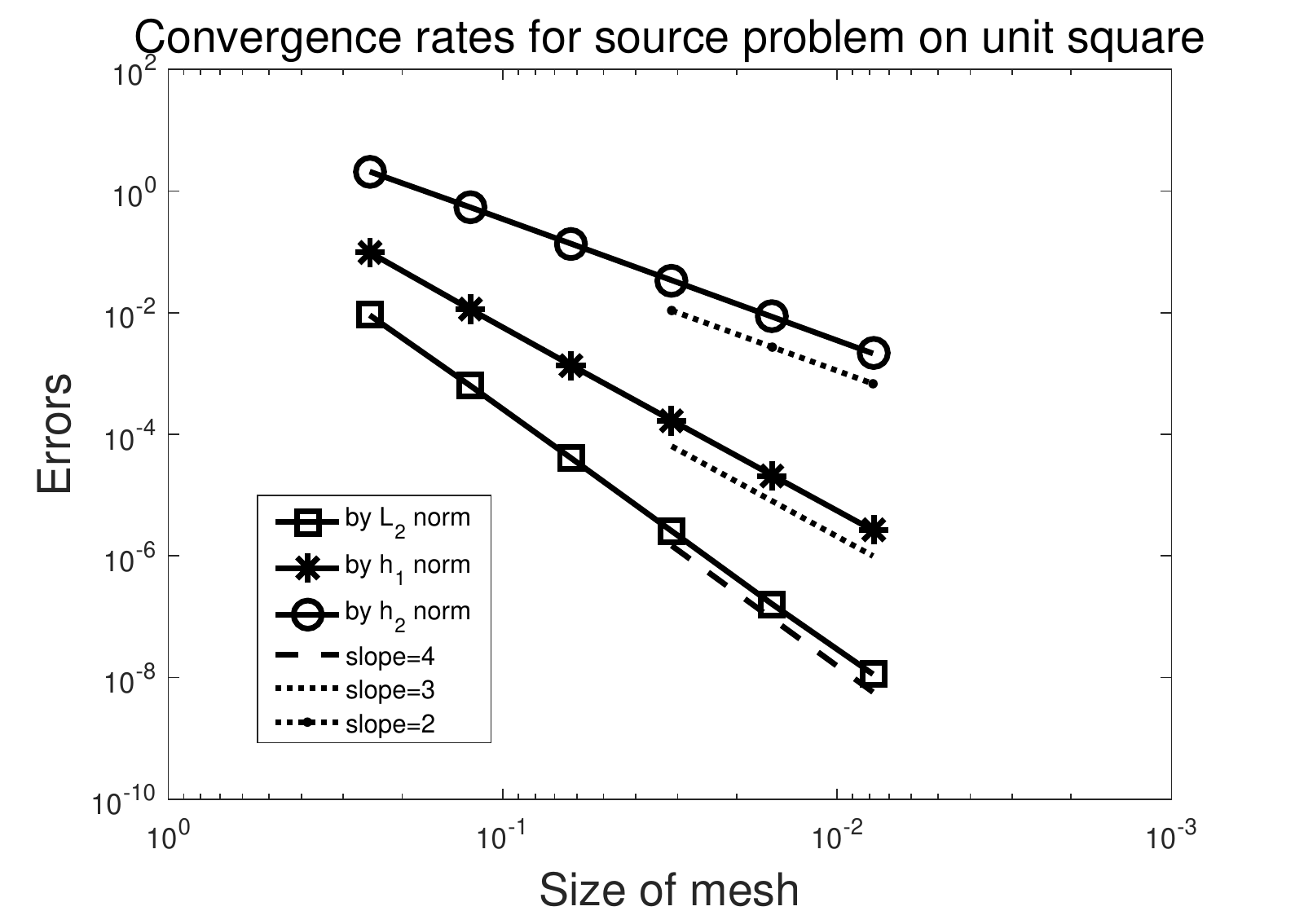}}
\subfigure[]{\label{SourceProblem_triangle}\includegraphics[scale=0.45]{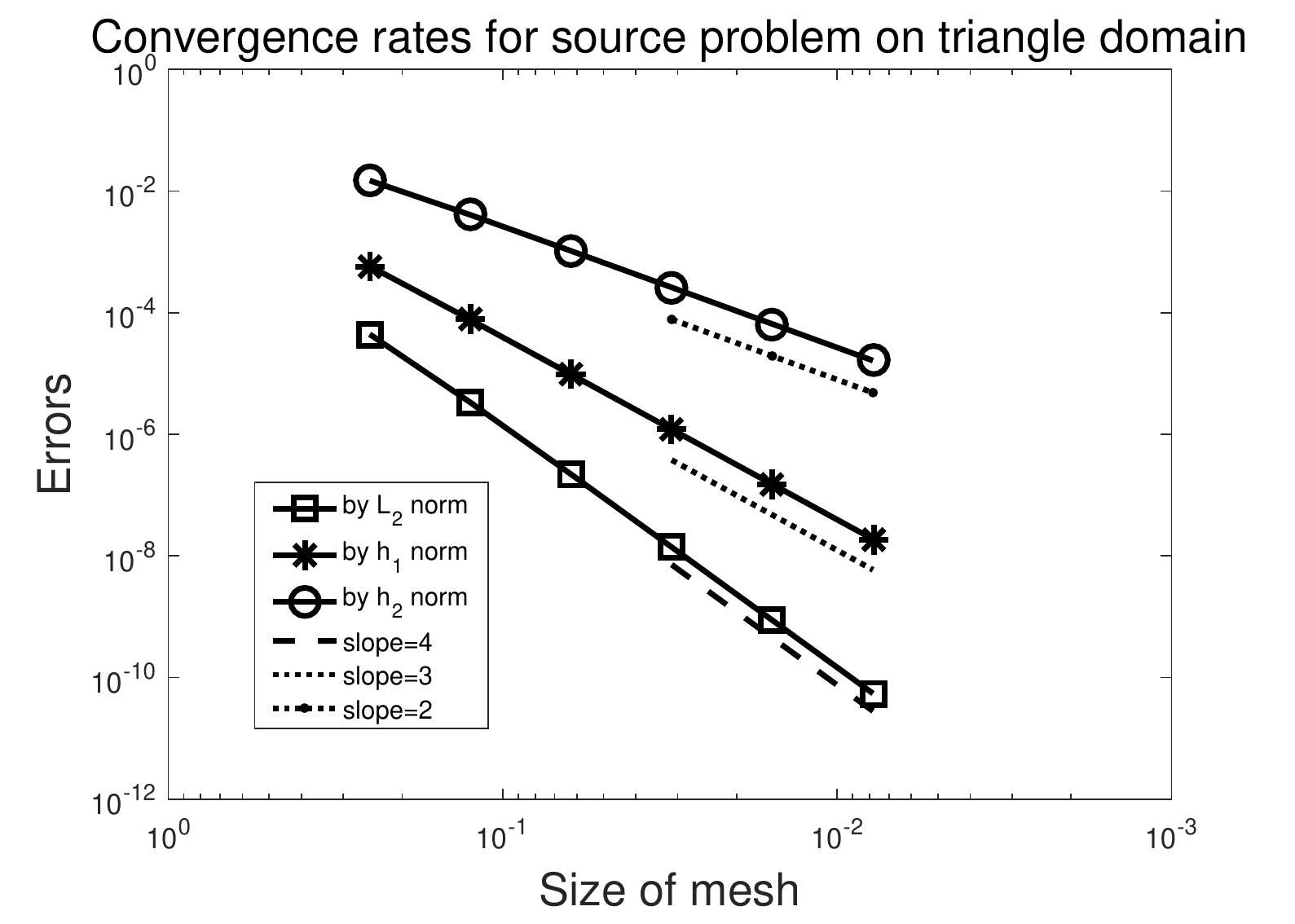}}
\caption{\label{fig:SourceProblem}The numerical performance for bi-Laplacian source problem by $B_{h0}^3$. Y-axis means the numerical error $\|u-u_{h_k}\|$ measured by $L_2$ or $h_1$ or $h_2$ norm. X-axis means the size of mesh. Left: for \textbf{Example 1} which is on square domain; Right: for \textbf{Example 2} which is on triangle domain.}
\end{figure}

For \textbf{Example 1}, the errors for numerical solutions are showed in Figure \ref{SourceProblem_square}. For \textbf{Example 2}, the errors for numerical solutions are showed in Figure \ref{SourceProblem_triangle}. We can observe that

(1) The convergence rate for source problem measured by $h_2$ norm is 2;

(2) The convergence rate for source problem measured by $h_1$ norm is 3;

(3) The convergence rate for source problem measured by $L_2$ norm is 4;

which are consistent with the theoretical results.

\textbf{Example 3.} Consider the bi-Laplacian source problem with varying coefficient $\delta=8+x_1-x_2$ on triangle domain $\Omega_2$ whose vertices are given by $(0,0), (1,0), (0,1)$.
And we consider $f=64x^3+48x^2y+528x^2-48xy^2+1152xy-368x-64y^3+624y^2-400y+64$ for which the exact solution is $u(x,y)=x^2y^2(1-x-y)^2$.

\begin{figure}[ht]
\begin{center}
\includegraphics[scale=0.75]{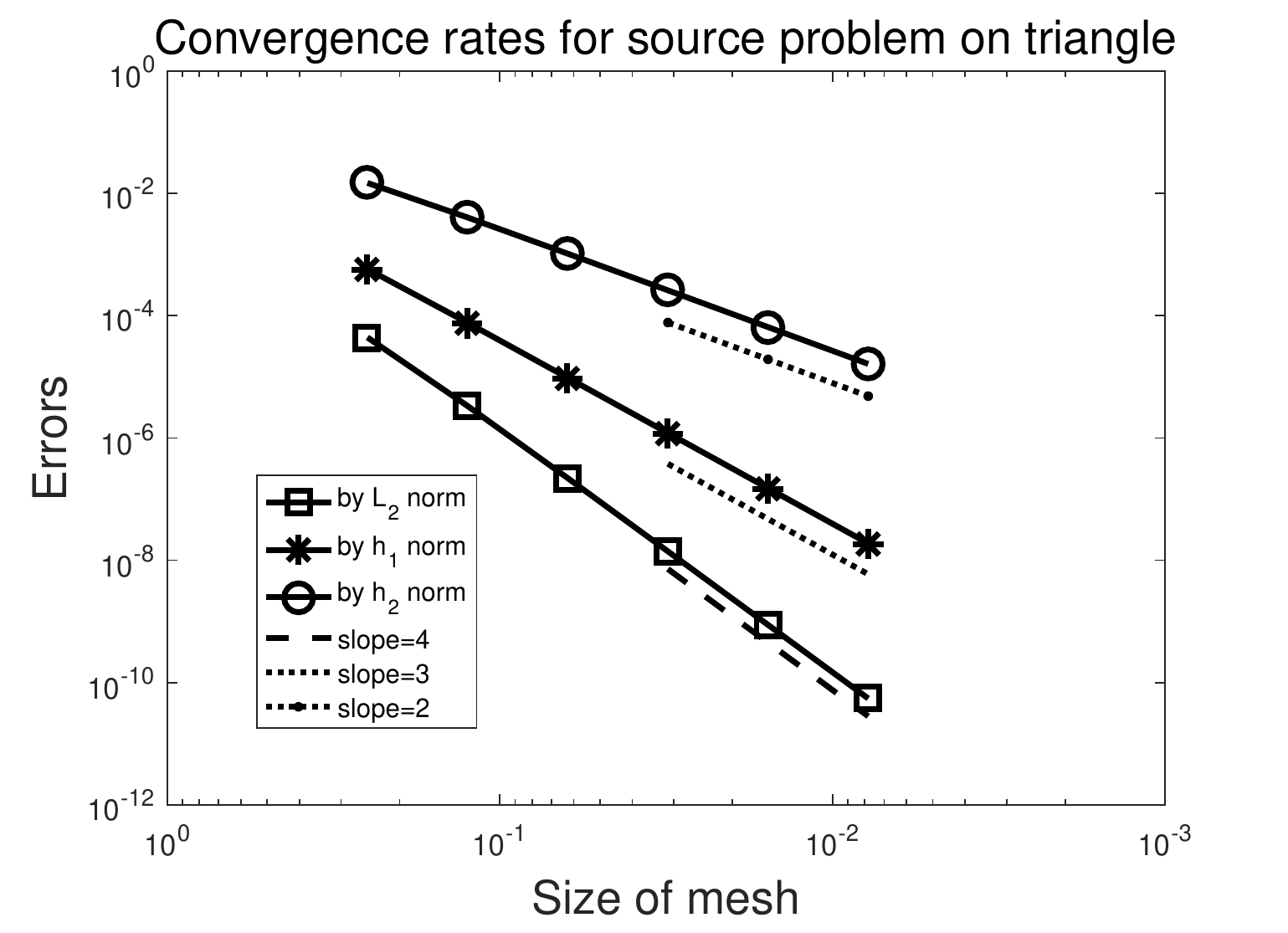}
\end{center}
\caption{The numerical performance by $B_{h0}^3$ for biharmonic source problem with non-constant coefficient $\delta=8+x_1-x_2$. Y-axis means the numerical error $\|u-u_{h_k}\|$ measured by $L_2$ or $h_1$ or $h_2$ norm. X-axis means the size of mesh.}
\label{BEP_VaryingCoefficient_triangle}
\end{figure}

For \textbf{Example 3}, the errors for numerical solutions are showed in Figure \ref{BEP_VaryingCoefficient_triangle}. We can observe that

(1) The convergence rate for source problem measured by $h_2$ norm is 2;

(2) The convergence rate for source problem measured by $h_1$ norm is 3;

(3) The convergence rate for source problem measured by $L_2$ norm is 4;

which are optimal and consistent with the theoretical results.

\subsubsection{For eigenvalue problem}

\textbf{Example 4.} Consider the bi-Laplacian eigenvalue problem (\ref{Non-constantBiharEigenProblem}) with constant coefficient $\delta=1$ on the unit square domain $\Omega_1=[0,1]^2$.

\textbf{Example 5.} Consider the bi-Laplacian eigenvalue problem (\ref{Non-constantBiharEigenProblem}) with constant coefficient $\delta=1$ on the non-convex L-shaped domain $\Omega_3=[0,1]\times[0,1]\backslash[0,\frac{1}{2})\times(\frac{1}{2},1]$.

\begin{figure}[ht]
\centering
\subfigure[]{\label{EigenProblem_square}\includegraphics[scale=0.505]{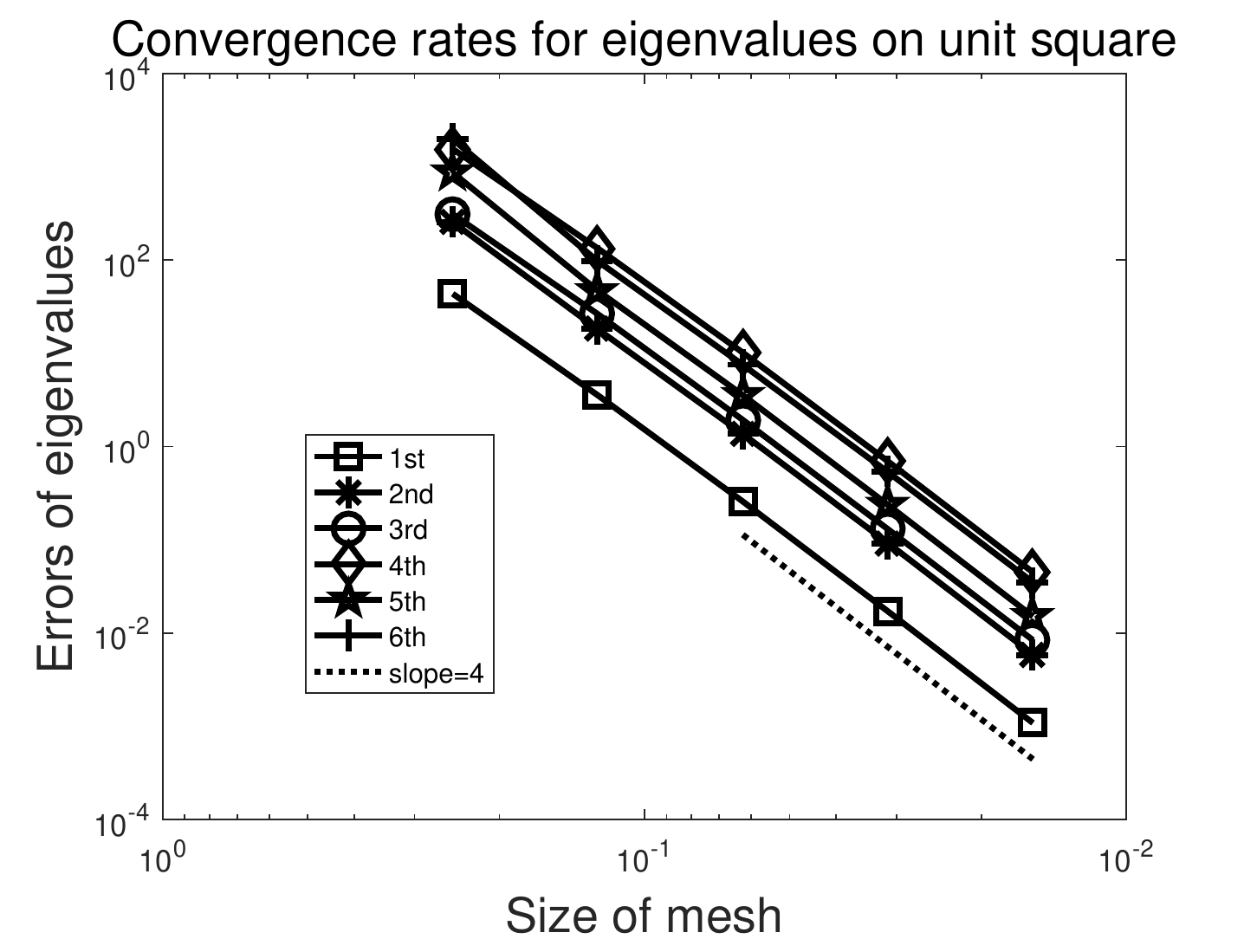}}
\subfigure[]{\label{EigenProblem_Lshaped}\includegraphics[scale=0.5]{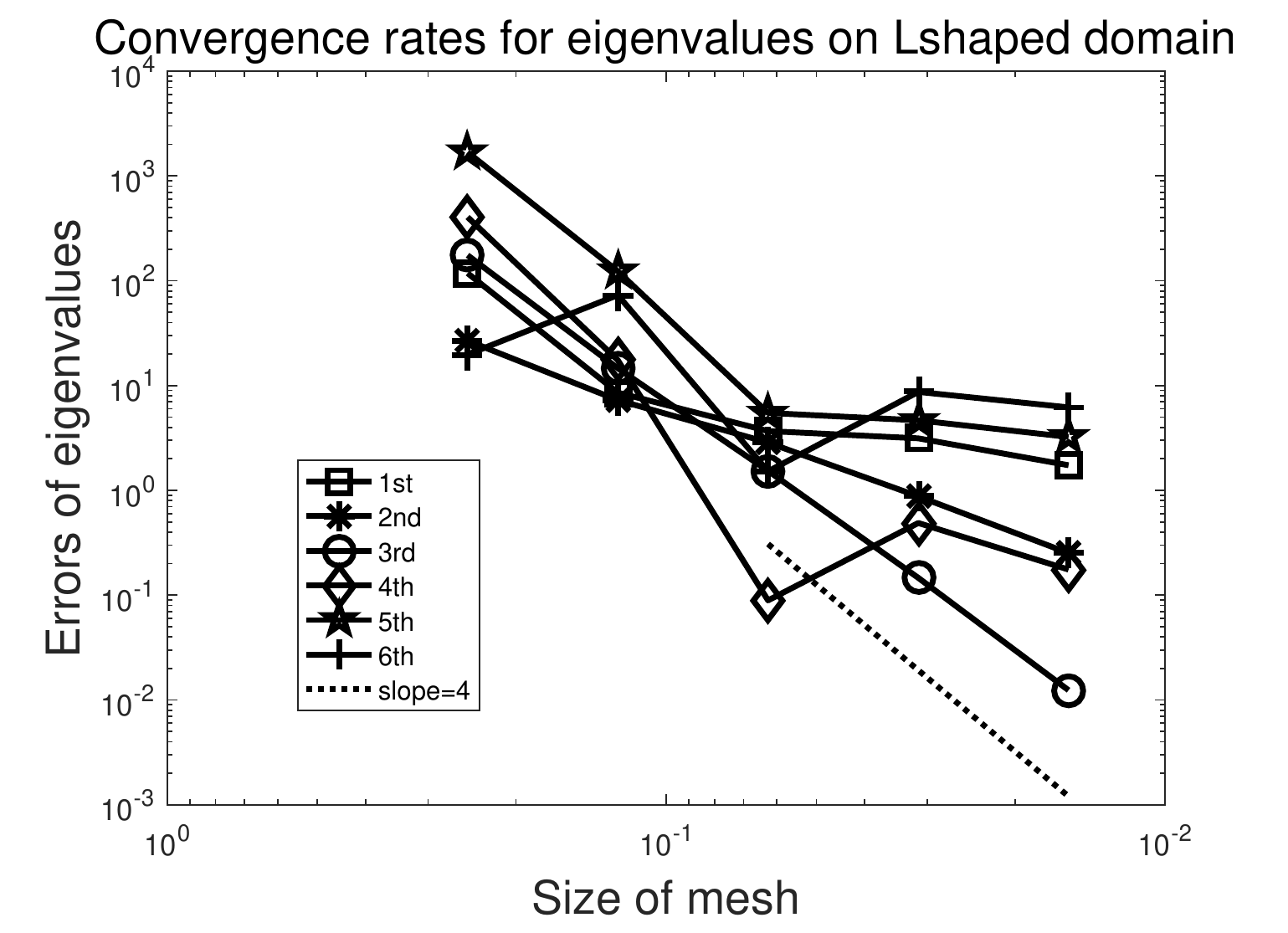}}
\caption{\label{BiharEigenProblemBh03}The convergence rates for the lowest six real eigenvalues for bi-Laplacian eigenvalue problem by $B_{h0}^3$. Y-axis means the numerical error $|\lambda-\lambda_{h_k}|$. X-axis means the size of mesh. Left: for \textbf{Example 4} which is on square domain; Right: for \textbf{Example 5} which is on the non-convex L-shaped domain.}
\end{figure}

Here we test on \textbf{Example 4} and \textbf{Example 5}, respectively. The mesh size of the
initial mesh is $h_0= \frac{1}{2}$. Six levels of uniformly refined triangular meshes are generated for numerical experiments and $h_k= h_{k-1}/2,\ k = 1,2,3,4,5,6$. The finest degrees of freedom (short for DOFs) for \textbf{Example 4} are 97283. The refinest DOFs for \textbf{Example 5} are 146435. We discretize by the second order computational scheme corresponding to $B_{h0}^3$ space. For each series of meshes, we obtain the computed eigenvalue $\lambda_{h_k}$. The convergent orders are computed by
\begin{equation}
log_2\left(|\frac{\lambda_{k-1}-\lambda_{k}}{\lambda_{k-2}-\lambda_{k-1}}|\right),\ \ \ \,k=3,4,5,6.\nonumber
\end{equation}
We present the results of the first six biharmonic eigenvalues showed in Figure \ref{BiharEigenProblemBh03}.
For \textbf{Example 4}, the results are showed in \ref{EigenProblem_square}. For \textbf{Example 5}, the numerical performance is showed in \ref{EigenProblem_Lshaped}. We can observe that for convex domain, the convergence rate for eigenvalues approximates 4 which is optimal and consistent with the theoretical expectation. For non-convex domain, the convergence rates are not optimal due to the low regularity of eigenfunctions.

\subsubsection{The $B_{h0}^3$ scheme for biharmonic eigenvalue problem with non-constant coefficient}
By $B_{h0}^3$ scheme, the variational formulation for (\ref{Non-constantBiharEigenProblem}) is as followed: find $u\in H_0^2(\Omega)$ and $\lambda\in R$, such that
\begin{equation}
(\delta\Delta u,\Delta v)=\lambda(u,v),\ \ \ \forall v\in H_0^2(\Omega),\nonumber
\end{equation}

The corresponding discretized variational formulation is to find $u_h\in B_{h0}^3$ and $\lambda_h\in R$, such that
\begin{equation}
(\delta\Delta u_h,\Delta v_h)=\lambda_h(u_h,v_h),\ \ \ \forall v_h\in B_{h0}^3.\nonumber
\end{equation}

\textbf{Example 6.} Consider the unit square domain $\Omega=[0,1]\times[0,1]$ with $\delta(x)=8+x_1-x_2$.

\textbf{Example 7.} Consider the unit square domain $\Omega=[0,1]\times[0,1]$ with $\delta(x)=\sqrt{x_1^2+x_2^2}+1$.

Here we test on \textbf{Example 6} and \textbf{Example 7}. The mesh size of the
initial mesh is $h_0= \frac{1}{4}$. Five levels of uniformly refined triangular meshes are generated for numerical experiments and $h_k= h_{k-1}/2,\ k = 1,2,3,4,5$. The finest degrees of freedom (short for DOFs) are 97283.

For \textbf{Example 6}, the lowest ten computed eigenvalues are showed in Table \ref{tab:example311}. The convergence rate of eigenvalues is 4. The computed eigenvalues tend to give the upper bound.

\begin{table}[htp]
\caption{\label{tab:example311}The performance of {\bf $B_{h0}^3$}  for \textbf{Example 6}.}
\begin{tabular}{cccccccc}\hline
Mesh&1&2&3&4&5&Trend&$Ord_{\lambda}$\\\hline
$\lambda_1$&10374.5195&10345.9954&10343.9256&10343.7882&10343.7794&$\searrow$&3.97049\\
$\lambda_2$&43152.3618&43005.8128&42994.7833&42994.0362&42993.9885&$\searrow$&3.96937\\
$\lambda_3$&43280.1536&43068.7439&43053.2500&43052.2064&43052.1391&$\searrow$&3.95288\\
$\lambda_4$&94720.7844&93650.3052&93568.7622&93563.1966&93562.8374&$\searrow$&3.95358\\
$\lambda_5$&138651.7814&138270.0014&138240.9393&138239.0035&138238.8805&$\searrow$&3.97531\\
$\lambda_6$&140390.6663&139603.9073&139543.0269&139538.7129&139538.4292&$\searrow$&3.92672\\
$\lambda_7$&221070.9885&217636.1630&217378.1490&217360.2947&217359.1410&$\searrow$&3.95190\\
$\lambda_8$&221623.7915&218016.9168&217724.6523&217704.3709&217703.0464&$\searrow$&3.93657\\
$\lambda_9$&353927.2977&353674.1751&353645.2752&353642.4924&353642.2935&$\searrow$&3.80689\\
$\lambda_{10}$&355323.7661&353783.9540&353664.7616&353656.2796&353655.7170&$\searrow$&3.91410\\\hline
\end{tabular}
\end{table}

For \textbf{Example 7}, the lowest ten computed eigenvalues are showed in Table \ref{tab:example312}. The convergence rate of eigenvalues is 4. The computed eigenvalues tend to give the upper bound.

\begin{table}[htp]
\caption{\label{tab:example312}The performance of {\bf $B_{h0}^3$}  for \textbf{Example 7}.}
\begin{tabular}{cccccccc}\hline
Mesh&1&2&3&4&5&Trend&$Ord_{\lambda}$\\\hline
$\lambda_1$&2242.0180&2236.1646&2235.7399&2235.7117&2235.7099&$\searrow$&3.97022\\
$\lambda_2$&9154.9841&9110.2819&9107.0020&9106.7807&9106.7664&$\searrow$&3.95159\\
$\lambda_3$&9486.4385&9456.1682&9453.8979&9453.7445&9453.7347&$\searrow$&3.97046\\
$\lambda_4$&20506.1886&20276.9808&20259.5801&20258.3940&20258.3174&$\searrow$&3.95238\\
$\lambda_5$&29816.2535&29732.0314&29725.6516&29725.2265&29725.1994&$\searrow$&3.97507\\
$\lambda_6$&30135.9791&29969.1169&29956.0394&29955.1090&29955.0478&$\searrow$&3.92537\\
$\lambda_7$&47066.3802&46285.9069&46222.8811&46218.5018&46218.2154&$\searrow$&3.93437\\
$\lambda_8$&49000.9675&48277.8380&48224.0164&48220.3162&48220.0777&$\searrow$&3.95558\\
$\lambda_9$&75673.2614&75431.9016&75400.3136&75398.0123&75397.8583&$\searrow$&3.90143\\
$\lambda_{10}$&75808.5004&75588.8154&75579.7274&75578.9302&75578.8748&$\searrow$&3.84894\\\hline
\end{tabular}
\end{table}

\subsection{Comparison with Morley element scheme}
\label{sec:morleybl}
We check the Morley element scheme for the eigenvalue problem
\begin{equation}
\left\{
\begin{array}{ll}
\Delta\delta\Delta u=\lambda u &\mbox{in}\,\Omega
\\
u=\frac{\partial u}{\partial n}=0,&\mbox{on}\,\partial\Omega.
\end{array}
\right.
\end{equation}
For Morley element, we consider the following variational formulation: find $u\in H_0^2(\Omega)$ and $\lambda\in R$, such that
\begin{equation}\label{NC_alphas_VariationalFormulation}
\alpha(\nabla^2 u,\nabla^2 v)+((\delta-\alpha) \Delta u,\Delta v)=\lambda(u,v),\ \ \ \forall v\in H_0^2(\Omega),
\end{equation}
where $(\nabla^2 u,\nabla^2 v)=\int_\Omega\sum_{s,t=1}^2\frac{\partial^2 u}{\partial x_s\partial x_t}\frac{\partial^2 v}{\partial x_s\partial x_t}dx,$ i.e., the inner product of the Hessian matrices of $u$ and $v$ and $\alpha$ is a constant satisfying $0<\alpha<\delta_{min}$. The items on the left side of (\ref{NC_alphas_VariationalFormulation}) guarantee the coercivity of variational problem.

The Morley element discretization space for $H_0^2(\Omega)$ is denoted by $V_h^M$. The corresponding discretized variational formulation is: find $u_h\in V_h^M$ and $\lambda_h\in R$, such that
\begin{equation}\label{Morley_alphas_VariationalFormulation}
\alpha(\nabla^2 u_h,\nabla^2 v_h)+((\delta-\alpha) \Delta u_h,\Delta v_h)=\lambda_h(u_h,v_h),\ \ \ \forall v_h\in V_h^M(\Omega).
\end{equation}

\begin{figure}[ht]
\begin{center}
\includegraphics[scale=0.75]{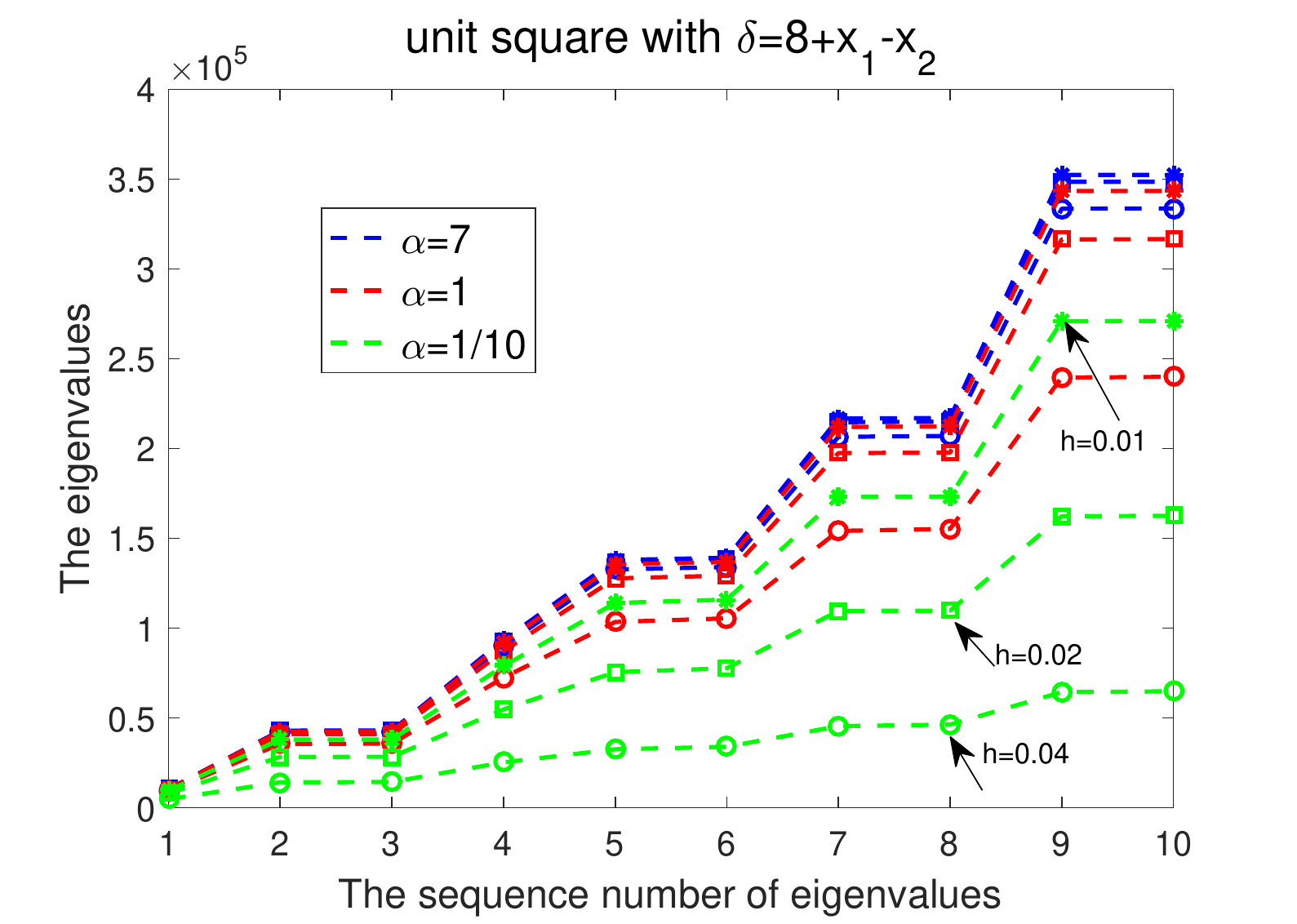}
\end{center}
\caption{The numerical performance by Morley element for biharmonic eigenvalue problem with non-constant coefficient $\delta=8+x_1-x_2$. Y-axis means the eigenvalues and X-axis means the sequence number of the lowest ten computed eigenvalues. For a fixed $\alpha$, the computed real eigenvalues on three successive grid levels are lised corresponding to mesh size $h=0.04, 0.02, 0.01$.}
\label{NBEP_alphas_nc1}
\end{figure}

We test the numerical performance of Morley element method on \textbf{Example 6} and \textbf{Example 7}.
For \textbf{Example 6}, by Morley element, the lowest ten computed real eigenvalues on three successive grid levels are showed in Figure \ref{NBEP_alphas_nc1}.
We can observe that the numerical results are sensitive to the parameter $\alpha$ and greatly depend on the choice of $\alpha$.

\begin{figure}[ht]
\begin{center}
\includegraphics[scale=0.75]{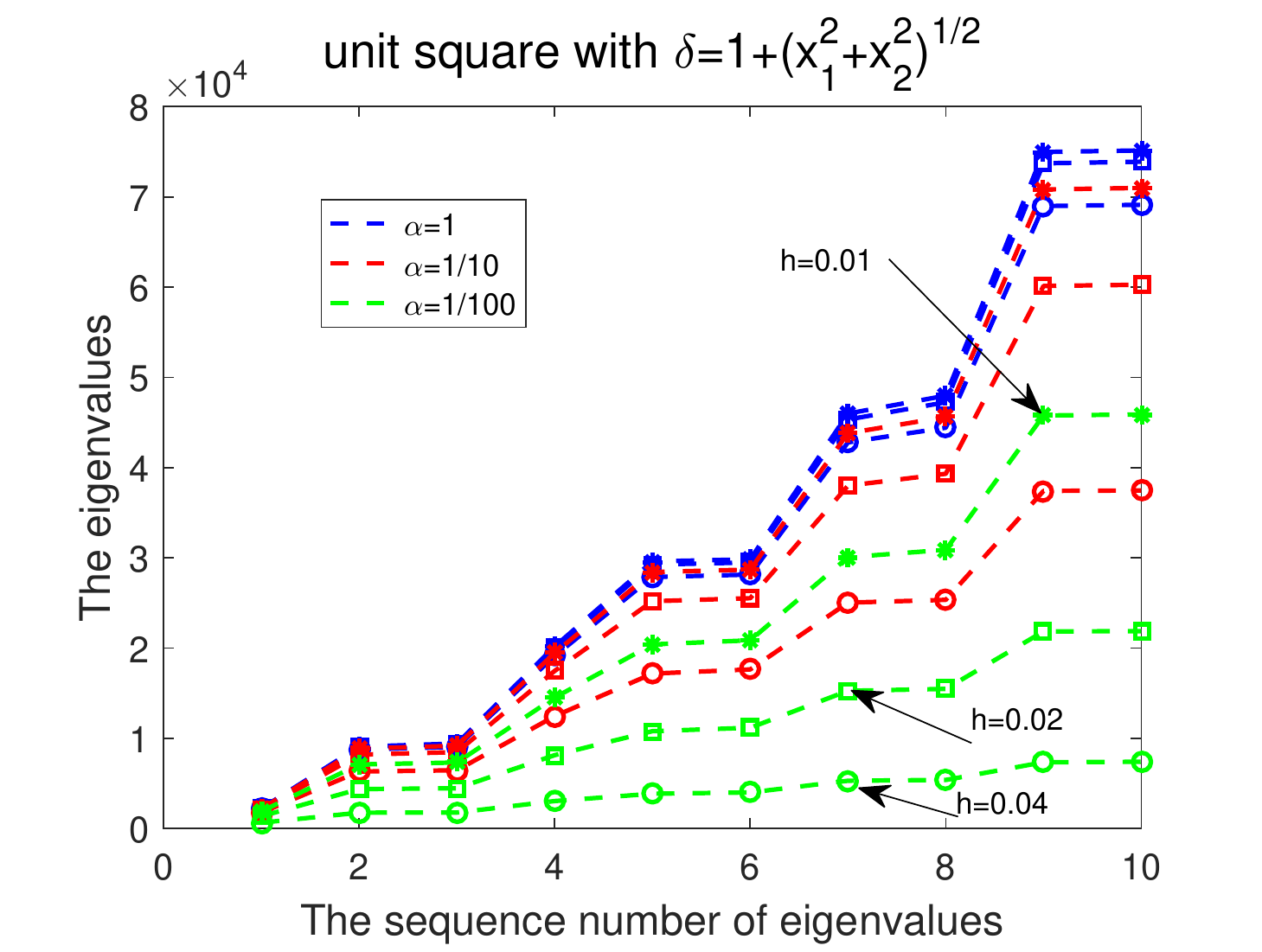}
\end{center}
\caption{The numerical performance by Morley element for biharmonic eigenvalue problem with non-constant coefficient $\delta=1+\sqrt{x_1^2+x_2^2}$. Y-axis means the eigenvalues and X-axis means the sequence number of the lowest ten computed eigenvalues. For a fixed $\alpha$, the computed real eigenvalues on three successive grid levels are lised corresponding to mesh size $h=0.04, 0.02, 0.01$.}
\label{NBEP_alphas_nc2}
\end{figure}

For \textbf{Example 7}, the numerical results are showed in Figure \ref{NBEP_alphas_nc2}. For different parameter $\alpha$, the computed eigenvalues are different. For different $\delta(x)$, the optimal $\alpha$ is also different.

\section{A high-accuracy scheme for the transmission eigenvalue problem}
\label{sec:te}


For the nonlinear transmission eigenvalue problem \eqref{Weak_Eigenvalue_A}, the corresponding discretization form is to find $(\tau_h,u_h)\in\mathcal{R}\times B_{h0}^3$ such that $\mathcal{B}(u_h,u_h)=1$ and
\begin{eqnarray}\label{Discrete_Eigenvalue_Problem_h}
\mathcal{A}_{\tau_h,h}(u_h,v_h) &=& \tau_h\mathcal{B}_h(u_h, v_h), \ \ \ \forall v_h\in B_{h0}^3.
\end{eqnarray}

Let $\{\xi_j\}^{N_h}_{j=1}$ be a basis for $B_{h0}^3$ and the corresponding FEM solution $u_h=\sum_{j=1}^{N_h} u_j\xi_j$, where $\{u_j\}$
corresponds to the standard degrees of
freedom for $B_{h0}^3$ scheme. We need the following matrices in the discrete case
\begin{table}[h]
\centering
\begin{tabular}{llc}
\hline
Matrix & Dimension & Definition\\ \hline  
$A$ & $N_h\times N_h$ &hessian matrix: $A_{i,j}=\int_\Omega\frac{1}{n-1}\Delta \xi_i\Delta \xi_jdx$\\
$B$ & $N_h\times N_h$ &stiff matrix: $B_{i,j}=\int_\Omega\frac{1}{n-1}\Delta\xi_i\xi_j+\frac{1}{n-1}\xi_i\Delta\xi_j-\nabla\xi_i\cdot\nabla\xi_j dx$ \\
$C$ & $N_h\times N_h$ & mass matrix: $C_{i,j}=\int_\Omega\frac{n}{n-1}\xi_i\xi_j dx$\\ \hline
\end{tabular}
\end{table}
and obtain the discretized quadratic eigenvalue problem
\begin{equation}\label{QEP}
(A+\tau B+\tau^2 C)x=0,
\end{equation}
where $x=(u_1,u_2,\cdots,u_{N_h})^T$. The computation of matrices $A,\ B,\ C$ involves the numerical integration of basis functions with non-constant coefficients. In practice, we use Gaussian integral formula and calculate the linear combination of function values at gaussian nodes on each triangular element.

For \eqref{QEP}, in practical computation, we convert to the linear eigenvalue problem
   \begin{equation}
     \left(
      \begin{array}{cc}
        -B & -A \\
         I &  O \\
      \end{array}
     \right)
     \left(
      \begin{array}{c}
      p_1  \\
      p_2  \\
      \end{array}
     \right)=\tau
     \left(
       \begin{array}{cc}
       C & O \\
       O & I \\
       \end{array}
     \right)
     \left(
       \begin{array}{c}
       p_1  \\
       p_2  \\
       \end{array}
     \right)\nonumber
   \end{equation}
and use matlab function "eigs" or "sptarn" to solve. And both $p_1$ and $p_2$ are all eigenvectors corresponding to $\tau$.

\begin{theorem}
Let $(\tau,u), (\tau_h,u_h)$ be the solution of (\ref{Weak_Eigenvalue_A}) and (\ref{Discrete_Eigenvalue_Problem_h}), respectively. Under the assumptions of Lemma 3.2 in \cite{Sun}, we can obtain the following results

\begin{eqnarray} \nonumber
\|u-u_h\|_{h_2}&\lesssim& h^2, \label{1}  \\ \nonumber
\|u-u_h\|_{1}&\lesssim& h^4, \label{2}   \\ \nonumber
|\tau-\tau_h|&\lesssim& h^4. \label{3} \nonumber
\end{eqnarray}
\end{theorem}

\subsection{The numerical performance of the nonconforming $B_{h0}^3$ scheme}
Here we focus on the case $n(x)>1$ which is of dominant interest in practice \cite{ColtonKress1998}. For the case $0<n(x)<1$, it can be treated similarly. Numerical experiments are conducted on a convex domain (a unit square domain $\Omega_1=[0,1]\times[0,1]$) and a non-convex domain (a L-shaped domain $\Omega_2=(0,1)\times(0,1)\backslash[\frac{1}{2},1)\times[\frac{1}{2},1)$). Six levels of uniformly refined triangular meshes are generated for numerical experiments. The mesh size of the initial mesh is $h_0=0.05$ and $h_k=h_{k-1}/2, k=1,2,3,4,5,6$. Note that further refinement would lead to very large matrix eigenvalue problems which take too long to solve. All examples are done using Matlab 2016a on a laptop with 16G memory and 2.9GHz Intel Core i7-7500U processor.

For each series of meshes, we obtain the eigenvalue series $\{\lambda_{h_k}\}_{k=1}^{6}$.
The convergent orders are computed by
\begin{equation}
log_2(|\frac{\lambda_{h_l}-\lambda_{h_{l+1}}}{\lambda_{h_{l+1}}-\lambda_{h_{l+2}}}|),\ \ \ \ l=1,2,3,4.
\end{equation}
We consider the following examples.

\textbf{Example 8.} The unit square domain $\Omega_1$ with the constant index of refraction $n(x)=16$.

The finest degrees of freedom (short for DOFs) are 194566. It costs 251.661752s for the whole calculation. We present the results of the first six real transmission eigenvalues. The eigenvalue approximations ($\lambda_h=\sqrt{\tau_h}$) on the finest mesh are (1.879591, 2.444236, 2.444236, 2.866439, 3.140111, 3.471509).

\begin{figure}[ht]
\begin{center}
\includegraphics[scale=0.55]{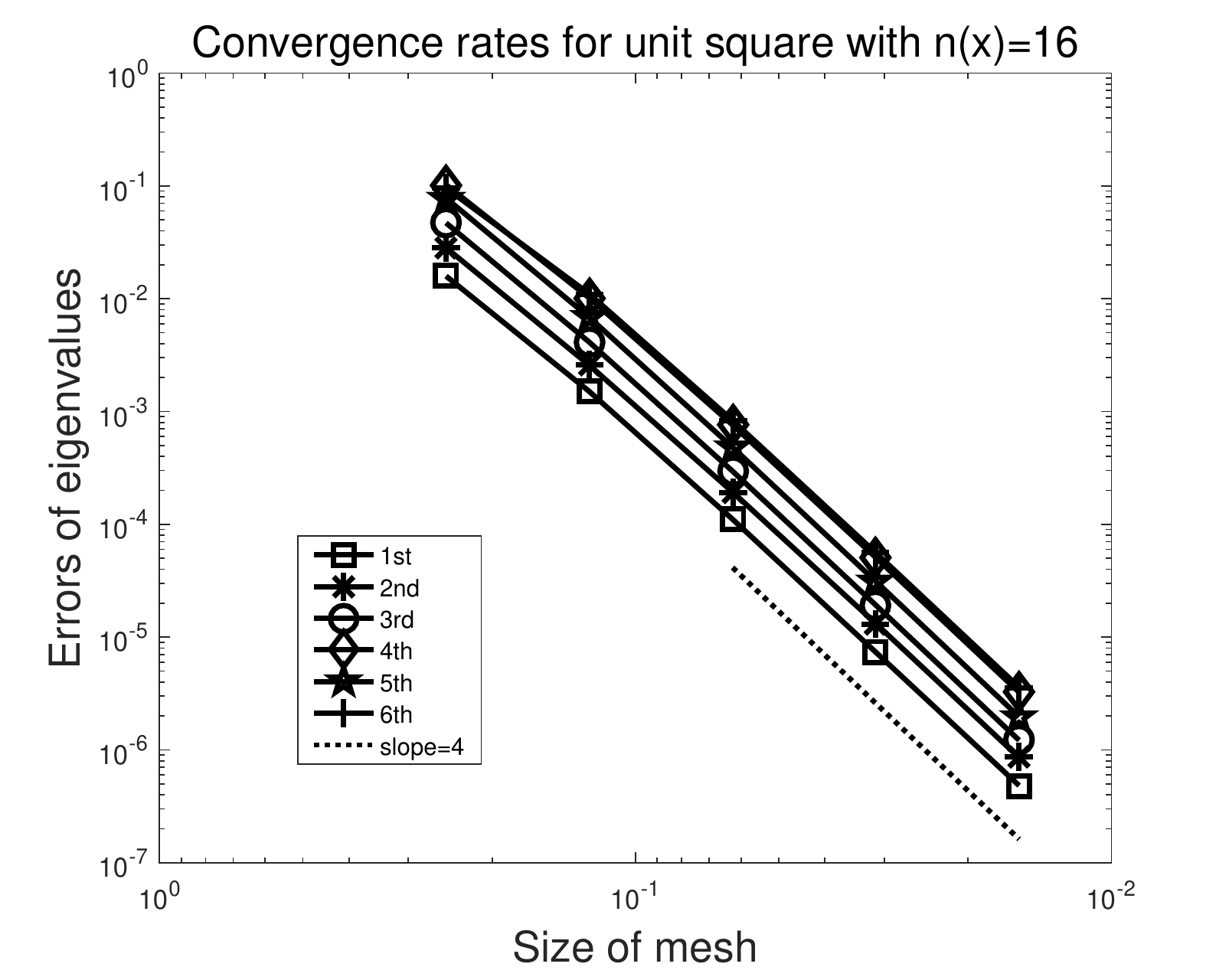}
\end{center}
\caption{The convergence rates for the lowest six real eigenvalues of the unit square with $n(x)=16$ by $B_{h0}^3$. Y-axis means $\lambda_{h_k}-\lambda_{h_6}$, as $h$ tends to zero, $\lambda_{h_k}-\lambda_{h_6}$ can be positive or negative; however, as illustrated in the figure, it's positive, here $k = 1,2,3,4,5$. X-axis means the size of mesh and so are the followings.}
\label{TEP_6evs_s16}
\end{figure}

From Figure \ref{TEP_6evs_s16}, we can observe the following phenomena:

(1) The convergence rates of transmission eigenvalues by $B_{h0}^3$ are 4.

(2) It gives the upper bound for real eigenvalues.

(3) The results by $B_{h0}^3$ are consistent with those in \cite{JiSunTurner2012ACMTOM}\cite{JiSunXie2014JSC}\cite{JiXiXie2017}.

\textbf{Example 9.} The unit square domain $\Omega_1$ with the non-constant index of refraction $n(x)=8+x_1-x_2$.

\begin{figure}[ht]
\begin{center}
\includegraphics[scale=0.55]{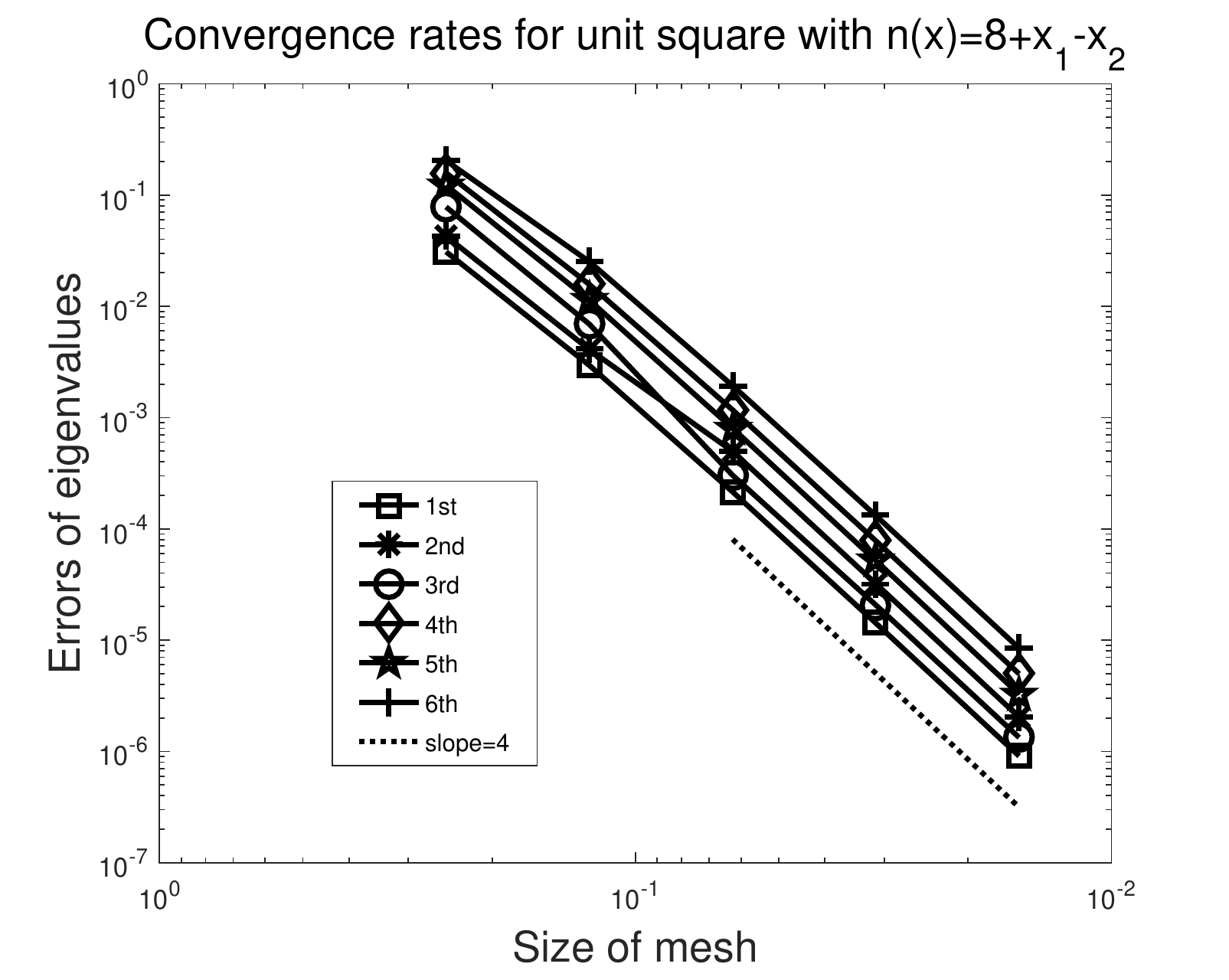}
\end{center}
\caption{The convergence rates for the lowest six real eigenvalues of the unit square with $n(x)=8+x_1-x_2$ by $B_{h0}^3$. Y-axis means $\lambda_{h_k}-\lambda_{h_6}$, as $h$ tends to zero, $\lambda_{h_k}-\lambda_{h_6}$ can be positive or negative; however, as illustrated in the figure, it's positive, here $k = 1,2,3,4,5$. X-axis means the size of mesh and so are the followings.}
\label{TEP_6evs_s_nc1}
\end{figure}

The first six real eigenvalue approximations on the finest mesh are (2.822189, 3.538697, 3.538992, 4.117742, 4.501729, 4.989140) which is consistent with the results in \cite{JiSunTurner2012ACMTOM}\cite{JiSunXie2014JSC}\cite{JiXiXie2017}. The convergence rates are showed in Figure \ref{TEP_6evs_s_nc1}. It can also be observed that $B_{h0}^3$ does give the theoretical predicted fourth convergence rate. And the computed real eigenvalues are monotonically decreasing as the mesh is refined.

\textbf{Example 10.} The L-shaped domain $\Omega_2$ with the constant index of refraction $n(x)=24$.

\begin{figure}[ht]
\begin{center}
\includegraphics[scale=0.55]{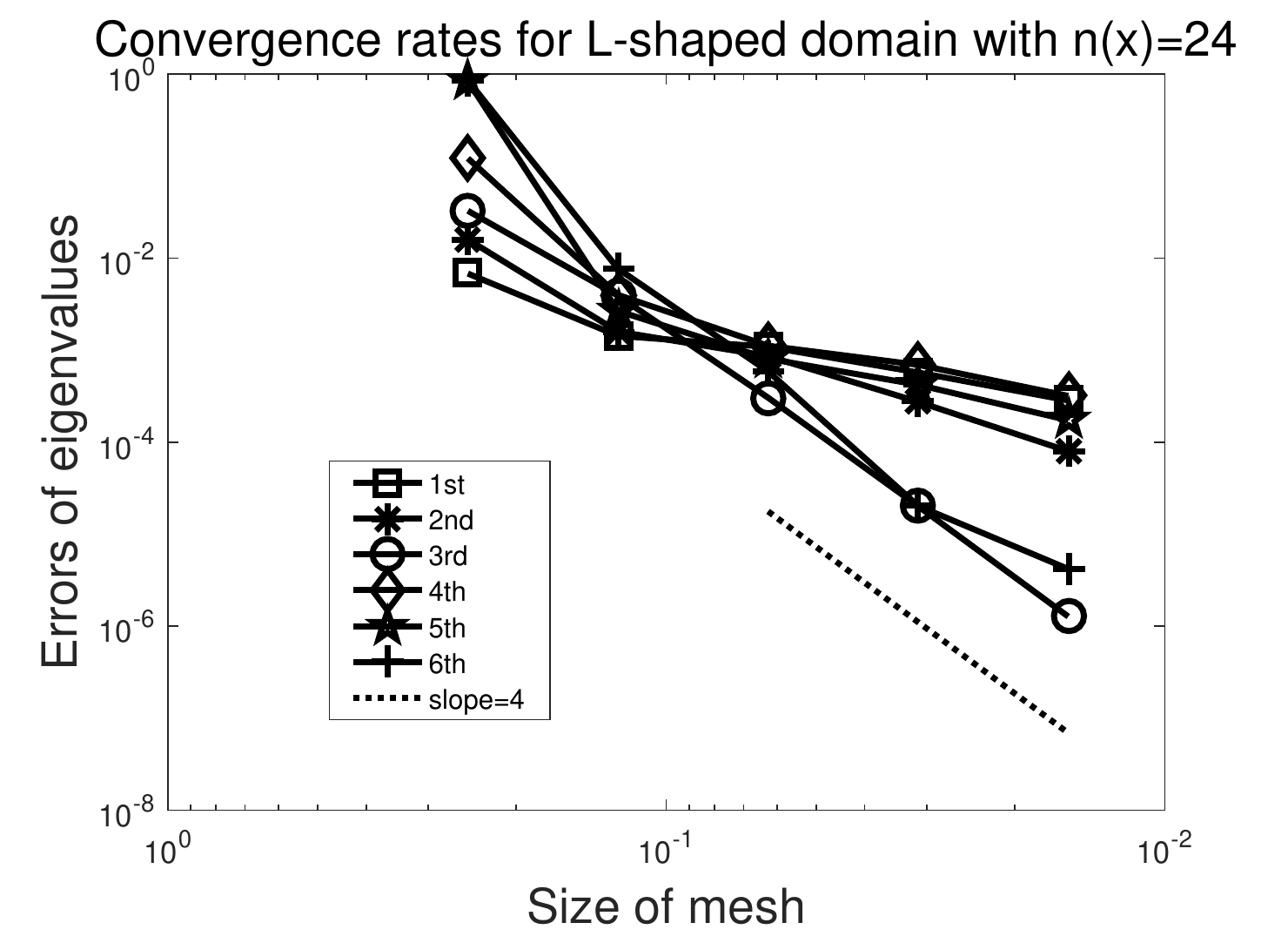}
\end{center}
\caption{The convergence rates for the lowest six real eigenvalues of the L-shaped domain with $n(x)=24$ by $B_{h0}^3$. Y-axis means $|\lambda_{h_k}-\lambda_{h_6}|$, as $h$ tends to zero. X-axis means the size of mesh.}
\label{TEP_6evs_L24}
\end{figure}

The finest DOFs are 292870. The total calculate time is 467.822844 second. The lowest six real eigenvalues on the finest mesh are (4.275620, 4.555635, 5.172225, 5.271284, 5.984808, 6.081556). Since $\Omega_2$ has a reentrant corner, the eigenfunction has a low regularity.
The convergence order for the eigenvalue approximation is less than 4 by the $B_{h0}^3$ scheme as is showed in Figure \ref{TEP_6evs_L24}.

\subsection{Morley element scheme revisited}
\label{sec:morleyte}

\begin{figure}[ht]
\begin{center}
\includegraphics[scale=0.75]{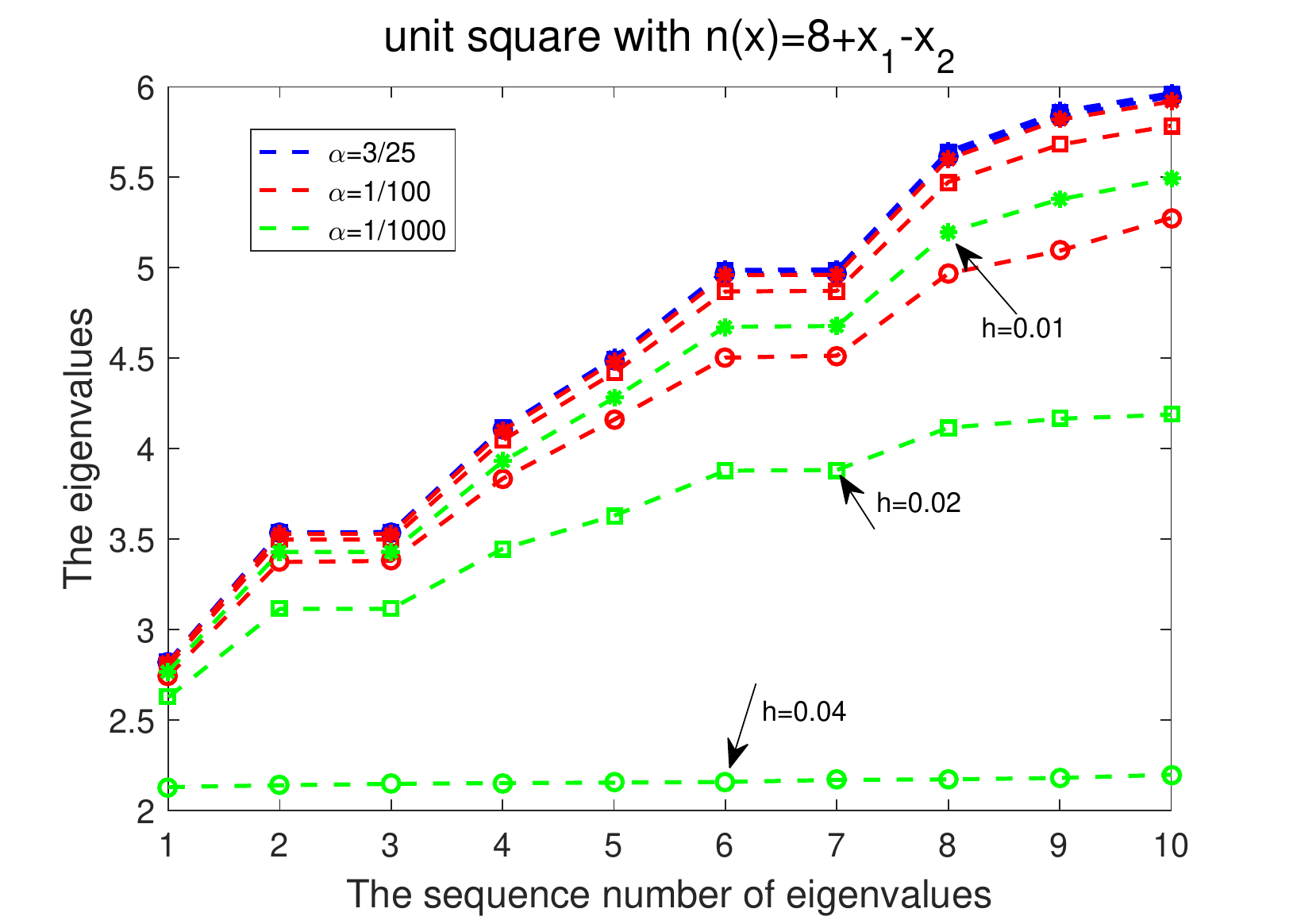}
\end{center}
\caption{The numerical performance by Morley element for transmission eigenvalue problem. Y-axis means the eigenvalues and X-axis means the sequence number of the lowest ten computed real eigenvalues. For a fixed $\alpha$, the computed real eigenvalues on three successive grid levels are lised corresponding to mesh size $h=0.04, 0.02, 0.01$.}
\label{Morley_alphas_nc1}
\end{figure}

\begin{figure}[ht]
\begin{center}
\includegraphics[scale=0.75]{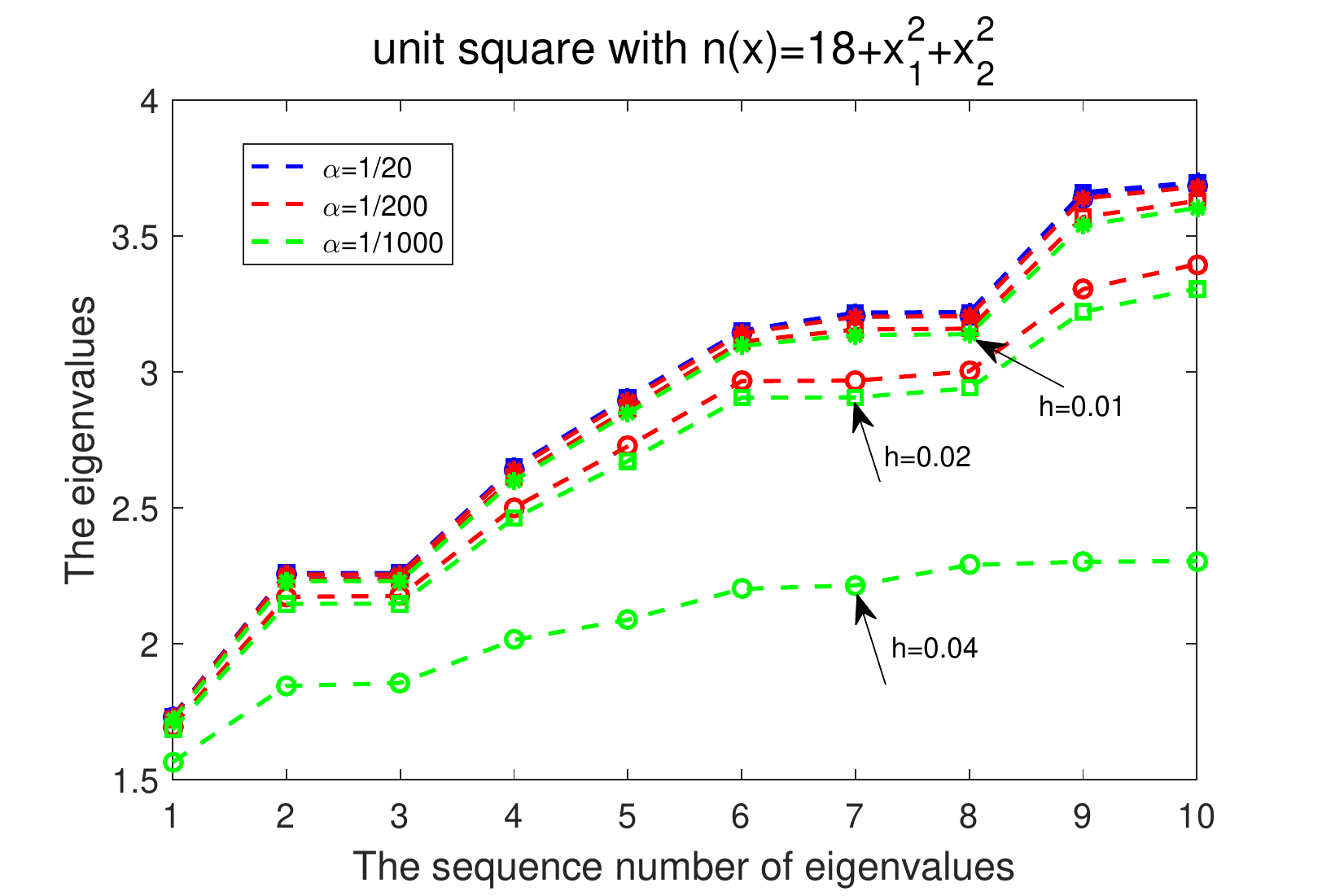}
\end{center}
\caption{The numerical performance by Morley element for transmission eigenvalue problem. Y-axis means the eigenvalues and X-axis means the sequence number of the lowest ten computed real eigenvalues. For a fixed $\alpha$, the computed real eigenvalues on three successive grid levels are lised corresponding to mesh size $h=0.04, 0.02, 0.01$.}
\label{Morley_alphas_nc2}
\end{figure}

In \cite{JiXiXie2017}\cite{XiJi2017}, the authors proposed the Morley element to discretize transmission eigenvalue problem. For the non-constant index of refraction $n(x)$, assume $0<\alpha_s\leq\frac{1}{n(x)-1}\leq\alpha_b$. They transformed the variational formulation to the following form:
\begin{equation}\label{revised_bilinear_form}
\Big(\frac{1}{n(x)-1}\Delta u,\Delta v\Big)=\Big(\Big(\frac{1}{n(x)-1}-\alpha\Big)\Delta u,\Delta v\Big)+(\alpha\nabla^2 u,\nabla^2 v),  
\end{equation}
where $(\nabla^2 u,\nabla^2 v)=\int_\Omega\sum_{s,t=1}^2\frac{\partial^2 u}{\partial x_s\partial x_t}\frac{\partial^2 v}{\partial x_s\partial
x_t}dx,$ i.e., the inner product of the Hessian matrices of $u$ and $v$, $\alpha$ is a constant satisfying $0<\alpha<\alpha_s$. The form on the right hand of (\ref{revised_bilinear_form}) guarantees the coercivity of the variational formulation on $H_0^2(\Omega)$ (c.f.\cite{XiJi2017}). However, in practical computation, the numerical performance is sensitive to the choice of $\alpha$. Figure \ref{Morley_alphas_nc1} shows the numerical performance by Morley element for unit square domain $\Omega=[0,1]^2$ with index of refraction
$n_1(x)=8+x_1-x_2$. We test on different $\alpha$. For a fixed $\alpha$, we record and present the lowest 10 computed real eigenvalues on three successive grid levels. It's observed that the numerical results are greatly dependent on the choice of $\alpha$.
Figure \ref{Morley_alphas_nc2} shows the numerical performance for unit square domain with index of refraction $n_2(x)=18+x_1^2+x_2^2$. For different index of refractions, the optimal choice of $\alpha$ is also different.

%
%
%
%
\section{Concluding remarks}
\label{sec:conc}

In this paper, we present a finite element scheme for the Helmholtz transmission eigenvalue problem based on the space $B^3_{h0}$. Different from most existing nonconforming finite elements, the bilinear form $(\Delta_h,\Delta_h)$ is coercive on the space $B^3_{h0}$, and it fits for the problem of operator $\Delta \delta\Delta$, including both the source and eigenvalue problems. Schemes associated with $B^3_{h0}$ are designed without introducing extra stabilisation mechanism. Numerical experiments illustrate the high accuracy of the schemes. Theoretical analysis will be given soon. The explicit formulation of the local basis functions obtained for easy application will bring in convenience in the future.


\begin{thebibliography}{99}
\bibitem{AnShen2013JSC} J. An and J. Shen,
	\textit{A {F}ourier-spectral-element method for transmission eigenvalue problems},
	J. Sci. Comput., \textbf{57} (2013), 670--688.
\bibitem{BabuskaOsborn} I. Babuska and J. E. Osborn,
 \textit{Eigenvalue Problems, in: P. G. Lions and P. G. Ciarlet (eds.) Handbook of Numerical Analysis, Vol. II, Finite Element Methods (Part 1)}, 641--787, NorthHolland, Amsterdam, 1991.
\bibitem{CakoniHaddar} F. Cakoni and H. Haddar, \textit{Transmission eigenvalues in inverse scattering theory}, Inside Out II, G. Uhlmann editor, MSRI Publications, \textbf{60} (2012), 526--578.
 \bibitem{CakoniRen} F. Cakoni and A. Ren, \textit{Transmission eigenvalues and the nondestructive testing of dielectrics}, Inverse Problems, 2008, 24, 065016 (15pp).
\bibitem{Camano2018} J. Camano, R. Rodriguez, P. Venegas, Convergence of a lowest-order finite element method for the transmission eigenvalue problem, Calcolo, \textbf{55} (2018).
\bibitem{ColtonKress1998} D. Colton and R. Kress,
    \textit{Inverse Acoustic and Electromagnetic Scattering Theory, 2nd ed.} Springer-Verlag, New York, 1998.
\bibitem{ColtonMonkSun2010} D. Colton, P. Monk, and J. Sun,
	\textit{Analytical and computational methods for transmission eigenvalues},
	Inver. Probl. \textbf{26} (2010), 045011.
\bibitem{Geng2016} H. Geng, X. Ji, J. Sun and L. Xu, \textit{C$^{0}$IP method for the transmission eigenvalue problem}, Journal of Scientific Computing, \textbf{68} (2016), 326-338.
\bibitem{HarrisCakoniSun} I. Harris, F. Cakoni and J. Sun,
 \textit{Transmission eigenvalues and non-destructive testing of anisotropic magnetic materials with voids}, Inverse Problems, 2014, 30: 035016 (21pp).
\bibitem{HuangEtal2016JCP} R. Huang, A. Struthers, J. Sun, and R. Zhang,
	\textit{Recursive integral method for transmission eigenvalues},
	J. Comput. Phys. \textbf{327} (2016), 830--840.
\bibitem{HuangEtal2017arXiv}R. Huang, J. Sun and C. Yang,
	\textit{Recursive integral method with Cayley transformation},
	arXiv:1705.01646, 2017.
\bibitem{JiLiu} X. Ji and H. Liu,
 \textit{On isotropic cloaking and interior transmission eigenvalue problems}, European Journal of Applied Mathematics, 2018, 29(2): 253-280.
\bibitem{JiSunTurner2012ACMTOM} X. Ji, J. Sun, and T. Turner,
	\textit{A mixed finite element method for Helmholtz transmission eigenvalues},
	 ACM Trans. Math. Software \textbf{8} (2012), Algorithm 922.
\bibitem{JiSunXie2014JSC}X. Ji, J. Sun, and H. Xie,
	\textit{ A multigrid method for Helmholtz transmission eigenvalue problem},
	J. Sci. Comput. \textbf{60} (2014), 276--294.
\bibitem{Kleefeld2013IP} A. Kleefeld,
	\textit{A numerical method to compute interior transmission eigenvalues}.
	Inver. Probl. \textbf{29} (2013), 104012.
\bibitem{JiXiXie2017} X. Ji, Y. Xi and H. Xie,
	\textit{Nonconforming Finite Element Method for the Transmission Eigenvalue Problem},
	Advances in Applied Mathematics and Mechanics, Vol.9, No.1, 2017, 92--103.
\bibitem{LiEtal2014JSC} T. Li, W. Huang, W. Lin, and J. Liu,
	\textit{On spectral analysis and a novel algorithm for transmission eigenvalue problems},
	J. Sci. Comput. \textbf{64} (2015), 83--108.	
\bibitem{Sun} J. Sun,
    \textit{Iterative methods for transmission eigenvalues},
    SIAM J. Numer. Anal., \textbf{49} (2011), 1860--1874.
\bibitem{XiJi2017} Y. Xi and X. Ji,
    \textit{Recursive integral method for the nonlinear non-selfadjoint transmission eigenvalue problem}, Journal of Computational Mathematics. Vol.35, No.6, 2017, 828--838.
    \bibitem{XiJi2018} Y. Xi, X. Ji, and H. Geng, A C$^0$IP Method of Transmission Eigenvalues for Elastic Waves, Journal of Computational Physics, \textbf{374} (2018), 237-248.
\bibitem{XiJiZhang2018} Y. Xi, X. Ji, and S. Zhang, A Multi-level Mixed Element Scheme for The Two
Dimensional Helmholtz Transmission Eigenvalue Problem, IMA Journal of Numerical
Analysis, accepted, 2018.
\bibitem{YangBiLiHan2016} Y. Yang, H. Bi, H. Li, and J. Han,
	\textit{Mixed methods for the Helmholtz transmission eigenvalues},
	SIAM J. Sci. Comput. \textbf{38} (2016), A1383--A1403.	
\bibitem{YangHanBi2016} Y. Yang, J. Han and H Bi,
    \textit{Non-conforming finite element methods for transmission eigenvalue problem}, Computer Methods in Applied Mechanics and Engineering, \textbf{307} (2016), 144-163.
\bibitem{Zhang2018} S. Zhang,
\textit{On optimal finite element schemes for biharmonic equation}, arXiv: 1805.03851.
\bibitem{Zhang2019} S. Zhang,
\textit{Optimal piecewise cubic finite element schemes for the biharmonic equation on general triangulations}, arXiv: 1903.04897.
\end{thebibliography}
\end{document}